\documentclass{amsart}

 0

\numberwithin{equation}{section}
\theoremstyle{plain}
	\begingroup
	\newtheorem{thm}[equation]{Theorem}
	\newtheorem{lem}[equation]{Lemma}
	\newtheorem{prop}[equation]{Proposition}

\endgroup
\theoremstyle{remark}
	\newtheorem{rem}[equation]{Remark}
\theoremstyle{definition}
	\begingroup
	\newtheorem{defn}[equation]{Definition}

	\newtheorem{algo}[equation]{}
	
\endgroup

\newsymbol\twoheadrightarrow 1310
\newcommand{\sobre}{\twoheadrightarrow}             
\newcommand{\inclu}{\hookrightarrow}
\newcommand{\ganda}{\rightharpoonup}                
\newcommand{\im}{\operatorname{Im}} 
\newcommand{\End}{\operatorname{End}}
\newcommand{\Aut}{\operatorname{Aut}}
\newcommand{\GL}{\operatorname{GL}}
\newcommand{\Ad}{\operatorname{Ad}}
\newcommand{\id}{\operatorname{id}}
\newcommand{\spn}{\operatorname{span}}
\newcommand{\Ind}{\operatorname{Ind}}
\newcommand{\bdng}{\operatorname{Braidings}}
\newcommand{\car}{\operatorname{char}}
\newcommand{\unid}[1]{{#1}^{\times}}
\newcommand{\prim}{\mathcal{P}}
\newcommand{\Ss}{\mathcal{S}}
\newcommand{\TTT}{{\mathbf T}}                          
\newcommand{\SSS}{{\mathbf S}}                          
\newcommand{\qqq}{{\mathbf q}}

\newcommand{\ZZ}{{\mathbb Z}}
\newcommand{\NN}{{\mathbb N}}
\newcommand{\BB}{{\mathbb B}}
\renewcommand{\SS}{{\mathbb S}}
\newcommand{\DD}{{\mathbb D}}
\renewcommand{\k}{{\mathbf k}}

\newcommand{\yd}{{}^H_H\mathcal{YD}}                
\newcommand{\ydskg}{{}^{\k\Gamma}_{\k\Gamma}\mathcal{YD}} 
\newcommand{\orb}{\mathcal{O}}                     
\newcommand{\com}[2]{#1_{(#2)}}                 
\newcommand{\nsq}[2]{(#1)_{\!#2}\,}             
\newcommand{\fsq}[2]{(#1)^{!}_{\!#2}\,}         
\newcommand{\csq}[3]{\left(\begin{array}{c}{#1}\\{#2}%
	\end{array}\right)_{\!\!#3}\,}
\newcommand{\bcp}[3]{\{#1_{#2},\ldots,#1_{#3}\}}
\newcommand{\vtn}[2]{#1^{\otimes #2}}           
\newcommand{\mdpd}[4]{\left(\begin{array}{rr}%
	#1 & #2 \\ #3 & #4\end{array}\right)}   
\newcommand{\mdpdc}[4]{\left(\begin{array}{cc}%
	#1 & #2 \\ #3 & #4\end{array}\right)}   
\newcommand{\trid}{\triangleright}                  

\newcommand{\Gm}{\Gamma}
\newcommand{\Dt}{\Delta}
\newcommand{\dt}{\delta}
\newcommand{\sgm}{\sigma}
\newcommand{\al}{\alpha}
\newcommand{\ep}{\epsilon}
\newcommand{\lmb}{\lambda}
\newcommand{\te}{\theta}
\newcommand{\ptl}{\partial}

\newcommand{\toba}[1]{{\mathfrak B}(#1)}
\newcommand{\tobag}[2]{{\mathfrak B}^{#1}(#2)}

\begin{document}

\headheight .08in
\title{On Nichols Algebras of low dimension}
\author{Mat\'\i as Gra\~na}
\address{Mat\'\i as Gra\~na \newline
\indent Depto. de Matem\'atica \newline
\indent Pab I, FCEyN, Ciudad Universitaria \newline
\indent (1428) Buenos Aires, Argentina}
\email{matiasg@dm.uba.ar}
\thanks{This work was partially supported by UBA, CONICOR and SeCyT-UNC}
\subjclass{Primary 16W30; Secondary 17B37}
\begin{abstract}
This is a contribution to the classification program of
pointed Hopf algebras. We give a generalization of the
quantum Serre relations and propose a generalization of
the Frobenius--Lusztig kernels in order to compute Nichols
algebras coming from the abelian case. With this,
we classify Nichols algebras $\toba V$ with dimension
$<32$ or with dimension $p^3$, $p$ a prime number,
when $V$ lies in a Yetter--Drinfeld category over
a finite group. With the so called Lifting Procedure, this
allows to classify pointed Hopf algebras of index $<32$ or
$p^3$.
\end{abstract}
\maketitle

\renewcommand{\theequation}{\arabic{section}.\arabic{subsection}.\arabic{equation}}
\section{Introduction, notation, preliminaries}
\subsection{Introduction}
The ``Lifting Procedure" (see \cite{as2}) is a method
intended to classify finite dimensional pointed Hopf
algebras. Among the achievements of this method, we
mention the classification of pointed Hopf algebras
of order $p^3,\ p^4$ ($p$ a prime number)
\cite{as2,as3}, the classification of pointed Hopf
algebras whose group is of exponent $p$ \cite{as4},
the classification of pointed Hopf algebras of order
$p^5$ \cite{as4,gp5,g32}, the classification
of pointed Hopf algebras of index $p,\ p^2$ \cite{gfr},
the answer, by the negative, to the tenth Kaplansky
conjecture \cite{as-cm} (some of these results have
been obtained independently by other methods, see for
instance \cite{bdg,ge,d,cd,cdr,svo,b}).
The first step in this procedure is the classification of
Nichols algebras, i.e., graded braided Hopf algebras with
the base field in degree $0$ and which are coradically
graded and generated by its primitive elements
(see \ref{df:nalg}). Being braided and graded, such an algebra
gives a braiding $c:V\otimes V\to V\otimes V$, where $V$
is the space of primitive elements (which turns out to
coincide with the homogeneous component of degree $1$).
Furthermore, the pair $(V,c)$ determines the Nichols
algebra up to isomorphism. A major goal in the classification
program of pointed Hopf algebras is thus to have an
algorithm describing (and in particular giving the
dimension of) a Nichols algebra starting out with the
pair $(V,c)$ as its space of primitive elements. Although
any solution to the Braid Equation $c:V\otimes V\to V\otimes V$
gives rise to a Nichols algebra, the objects $V$ arising
in the application of the Lifting Procedure are objects
of Yetter--Drinfeld categories over finite groups, and
hence we shall restrict ourselves to those pairs $(V,c)$
with this additional property (see \ref{df:pttg}).

The theory divides at this moment in two branches: the pair
$(V,c)$ is ``of diagonal group type" (see \ref{df:pttg}) or
it is not. In the first case, a substantial progress is made
in \cite{as-cm} through the notion of ``braiding of Cartan
type", using the theory of Frobenius--Lusztig kernels.
This notion applies to a large class of pairs $(V,c)$
of diagonal group type. In this article we propose to
enlarge the class of Frobenius--Lusztig kernels and we
give the examples of type $A_2$ and $B_2$ (see \ref{pr:mdr2}
and \ref{rm:cpnc}).  To compute the PBW bases of these
generalized FL-kernels we use a generalized form of the
quantum Serre relations.  Though this generalized form
can be found in \cite{ros} for the diagonal case, we
prefer to derive it from an equality in the group algebra
of the braid group, since this equality could provide
quantum Serre relations for the general case (see
\ref{pr:qsrg} and \ref{th:qsrg}).

In the non-diagonal case very little is known. In fact,
up to now there were only two examples known to produce
finite dimensional Nichols algebras (and hence, finite
dimensional pointed Hopf algebras). The first one is
$3$-dimensional and it generates a $12$-dimensional
Nichols (it is given for the first time in
\cite{ms} and was reproduced in \cite{ag}). The second
one is $6$ dimensional and it is given also in
\cite{ms}. Another example given there, which is $4$
dimensional, is proved in this article to be of diagonal
group type (see \ref{rm:msvca}).\footnote{After submitting
the paper I was informed that a third module given in
\cite{ms} is known to produce a finite dimensional Nichols
algebra (see \cite{fk}).}
In this article we propose some tools to classify all
Yetter--Drinfeld modules of a given dimension over
finite groups. We use these tools to parameterize those
modules of dimension $3$, $4$ (sections \ref{ss:r3},
\ref{ss:r4}). As a consequence of this, we find a new
example of a finite dimensional Nichols algebra with a
$4$-dimensional space of primitive elements which is
of non-diagonal group type; its dimension is $72$
(see \ref{al:cdt72}).
We also give some necessary conditions for an irreducible
YD-module to be of ``finite Nichols rank" (i.e., to give a
finite dimensional Nichols algebra).
This allows to focus the attention mainly on the specific
cases $V=M(g,\chi)$ where $\chi$ is a character (see
\ref{lm:csgr} for the condition, \ref{df:mgr} for the
definition of $M(g,\rho)$).

Finally, as a consequence of these methods together with
those of \cite{gfr}, we classify Nichols algebras over
finite groups with dimension $<32$. As said above, this
is the first step in the classification of pointed Hopf
algebras of index $<32$. Specifically, we have:
\begin{thm}
Let $R\in\ydskg$ be a Nichols algebra, where $\Gm$ is a
finite group and $\dim R<32$. Let $V=\prim R$ be its space
of primitive elements. Then $\dim V\le 4$ and either
\begin{itemize}
\item $R$ is a Quantum Linear Space (see definition in \ref{df:qls}), or
\item $R$ is of type $A_2$ (see remark \ref{rm:cpnc}), or
\item $R$ is of type $A_2\times A_1$, or
\item $V=V_{3,-1}$ and $\dim R=12$ (see \ref{df:v3q} for the definition), or
\item $V=V_{3,-1}\oplus V'$, $\dim\toba{V'}=2$,
        $c^2:V_{3,-1}\otimes V'\to V_{3,-1}\otimes V'$ is the identity,
        and $\dim R=24$.
\end{itemize}
\end{thm}

\subsection{Acknowledgements}
I thank N. Andruskiewitsch for the careful reading of the
material and for general and particular advice. I thank
J. Fern\'andez Bonder for a hint on \ref{rm:osf}. The referee
has suggested the best terminology I have seen for the braided
pairs that this article deals with. It has been a pleasure
to read his comments.
I thank also M. Lamas for many valuable conversations.

\subsection{Notation}\label{ss:bp}
For a general reference on Hopf algebras see \cite{mon}.
We fix an algebraically closed field $\k$ of characteristic
$0$. Most of the results of the article are still true for
positive characteristic, as long as it does not divide the
dimension of the Nichols algebras under consideration (see
\cite{gtesis}). In this case, the definition of $N(1)$ below must
be given as $N(1)=\car(\k)$.

All vector spaces, algebras, tensor products, homs are
considered over $\k$ unless explicitly stated. If $\Gm$ is
a group and $g\in\Gm$, we denote by $\Gm_g$ the centralizer
$\Gm_g=\{x\in\Gm\ |\ xg=gx\}$.

For a root of unity $q$, we denote by $N(q)$ its order if
$q\neq 1$. If $q=1$ or $q$ is not a root of unity then we
define $N(q)=\infty$. We recall the definition of $q$-numbers:
for $\qqq$ an indeterminate and $n,m\in\NN$, $n\ge m$, define
$$\nsq n{\qqq}=1+\qqq+\qqq^2+\cdots+\qqq^{n-1},\quad
\fsq n{\qqq}=\prod_{i=1}^n\nsq i{\qqq},\quad
\csq nm{\qqq}=\frac{\fsq n{\qqq}}
	{\fsq m{\qqq}\fsq {n-m}{\qqq}}.$$
All these expressions lie in $\ZZ[\qqq]$. For $q\in\k$ we
evaluate them in $\qqq=q$ to get $\nsq nq$, $\fsq nq$ and
$\csq nmq$.

Let $\SS_n$ be the symmetric group, presented with $n-1$
generators $\{\tau_1,\ldots,\tau_{n-1}\}$ and relations
$\tau_i\tau_j=\tau_j\tau_i$ if $|i-j|>1$,
$\tau_i\tau_{i+1}\tau_i=\tau_{i+1}\tau_i\tau_{i+1}$,
$\tau_i^2=1$. We denote by $\BB_n$ the braid group, with
$n-1$ generators $\{\sgm_1,\ldots,\sgm_{n-1}\}$ and
relations
$$\sgm_i\sgm_j=\sgm_j\sgm_i\mbox{ if }|i-j|>1,\quad
	\sgm_i\sgm_{i+1}\sgm_i=\sgm_{i+1}\sgm_i\sgm_{i+1}.$$
We consider the sequence of inclusions
$$\BB_1\inclu \BB_2\inclu\ldots\inclu\BB_n
	\inclu\BB_{n+1}\inclu\ldots$$
where $\sgm_i\mapsto\sgm_i$, and consider $\BB_{\infty}$
the limit of this sequence.

\begin{defn}
For $m<n$, we denote by $\BB_{m\_n}$ the subgroup of
$\BB_n$ generated by $\sgm_m,\ldots,\sgm_{n-1}$. We denote by
$\iota_{m\_n}:\BB_{n-m+1}\to\BB_n$ the inclusion $\sgm_i\mapsto\sgm_{i+m-1}$.
Then $\BB_{m\_n}=\im(\iota_{m\_n})$. We may view all these
groups as subgroups of $\BB_{\infty}$. If $x_1,x_2$ are elements
of $\BB_{j_1},\BB_{j_2}$ respectively, then we denote
$$(x_1|x_2)=x_1\cdot\iota_{(j_1+1)\_(j_1+j_2)}(x_2)\in\BB_{j_1+j_2},$$
and inductively we define
$$(x_1|\cdots|x_s)=(x_1|(x_2|\cdots|x_s)) \in\BB_{j_1+\cdots+j_s}.$$
For instance, if $x_1=\sgm_2\sgm_1\in\BB_3$ and $x_2=\sgm_1\in\BB_2$,
we have $(x_1|x_2)=\sgm_2\sgm_1\sgm_4\in\BB_5$.
We extend (bi)linearly the definition of $(|)$ to the group algebra,
i.e.
$$(|):\k\BB_{j_1}\otimes\k\BB_{j_2}\to\k\BB_{j_1+j_2}.$$
\end{defn}

\subsection{Preliminaries}
We consider the canonical projection $\BB_n\sobre\SS_n$,
$\sgm_i\mapsto\tau_i$. By \cite[64.20]{cure} we have a
section (of sets) $s:\SS_n\to\BB_n$ such that if
$x=\tau_{i_1}\cdots\tau_{i_l}$, where $l=\ell(x)$ (i.e.,
$x$ can not be written with fewer than $l$ generators
$\tau_i$'s), then $s(x)=\sgm_{i_1}\cdots\sgm_{i_l}$.
Notice that if $\ell(xy)=\ell(x)+\ell(y)$ then
$s(xy)=s(x)s(y)$.

We denote by $\SSS^n$ the element of the group algebra $\k\BB_n$
\begin{equation}\label{df:ddsn}
\SSS^n=\sum_{x\in\SS_n}s(x).
\end{equation}
Let $i_1+\cdots+i_r=n$. We denote by $Sh_{i_1,\ldots,i_r}\subseteq\SS_n$
the subset of $(i_1,\ldots,i_r)$-shuffles, i.e.
\begin{multline*}
Sh_{i_1,\ldots,i_r}=\{x\in\SS_n\ |\ 
	x^{-1}(i_1+\cdots+i_j+1)<x^{-1}(i_1+\cdots+i_j+2)<\cdots \\
	\cdots <x^{-1}(i_1+\cdots+i_j+i_{j+1})\mbox{ for }j=0,\ldots,r-1\}.
\end{multline*}
We consider also the elements
\begin{align*}
\SSS_{i_1,\ldots,i_r} &= \sum_{x\in Sh_{i_1,\ldots,i_r}}s(x)
	\in\k\BB_{i_1+\cdots+i_r}, \\
\TTT_{i_1,\ldots,i_r} &= \sum_{x^{-1}\in Sh_{i_1,\ldots,i_r}}s(x)
	\in\k\BB_{i_1+\cdots+i_r}, \\
\SSS^{j_1|j_2|\cdots|j_s} &= (\SSS^{j_1}|\SSS^{j_2}|\cdots|\SSS^{j_s})
	\in\k\BB_{j_1+\cdots+j_s}.
\end{align*}
Notice that
$$\SSS^{i|j}\SSS_{i,j}=\SSS_{i+j}=\TTT_{i,j}\SSS^{i|j}.$$

\medskip
\begin{defn}
We call a {\it braided pair} (or BP) a pair $(V,c)$, where
$V$ is a vector space and $c:V\otimes V\to V\otimes V$ is
a (bijective) solution of the braid equation. We say that
$(V,c)$ is a {\it rigid braided pair} (or RBP) if it is a
finite dimensional BP such that
$c^{\flat}\in\End(V^*\otimes V,V\otimes V^*)$ is an isomorphism
(see for instance \cite{tak1} for the definition of $c^{\flat}$).
\end{defn}

Let $H$ be a Hopf algebra with bijective antipode.
We denote by $\yd$ the category of (left-left)
Yetter--Drinfeld modules over $H$ of arbitrary dimension.
If $V$ is an object of $\yd$ with braidings
$c_{V,W}:V\otimes W\to W\otimes V$, then $(V,c=c_{V,V})$ is a braided pair.
If $V$ is moreover finite dimensional then $(V,c=c_{V,V})$ is an RBP.

Let $(V,c)$ be a BP. The group $\BB_n$ acts on $\vtn Vn$ by
$\sgm_i\mapsto(\id_{\vtn Vi}\otimes c\otimes\id_{\vtn V{n-i-2}})$.
We denote respectively by $AV,CV$ the tensor algebra
and the tensor coalgebra of $V$. As vector spaces, they coincide with
$TV=\oplus_i\vtn Vi=\oplus_iA^iV=\oplus_iC^iV$.\linebreak
The $(i,j)$-component of the multiplication of $AV$ is simply the
identity map \linebreak
$m_{i,j}=\id:A^iV\otimes A^jV\to A^{i+j}V$.
Dually, the $(i,j)$-component of the comultiplication of $CV$ are
the identity maps $\Dt_{i,j}=\id:C^{i+j}V\to C^iV\otimes C^jV$.
Both $AV$ and $CV$ are graded braided Hopf algebras.
The comultiplication components of $AV$
are nothing but $\SSS_{i,j}:A^{i+j}V\to A^iV\otimes A^jV$.
Dually, the multiplication components of $CV$ are
$\TTT_{i,j}:C^iV\otimes C^jV\to C^{i+j}V$. The antipode in both cases
is determined by $\Ss(x)=-x\ \forall x\in V$.

There exists a unique (graded) map of braided Hopf algebras which reduces
to the identity $\id:V\to V$ in degree $1$; we denote this map by $\SSS$.
Its graded components are the maps $\SSS^n:A^nV\to C^nV$.

\begin{defn}\label{df:nalg}
Let $(V,c)$ be an RBP. We say that the image of
$$\SSS:AV\to CV$$
is the {\em Nichols algebra} generated by $V$ and we denote it by
$\toba V$ \footnote{If $(V,c)$ is just a BP, the object $\toba V$
can be constructed in the same manner. In this general case, the
result has been called a {\em Quantum Symmetric Algebra}.}.
We say that $\dim(\toba V)$ is the {\em Nichols rank} of $V$.
\end{defn}

\begin{rem}
It can be proved that a finite dimensional braided pair $(V,c)$
belongs to a braided category, and that an RBP $(V,c)$ belongs
to a rigid braided category. We warn that in \cite{as-cm} the
term ``Nichols algebra" is used for the case in which the rigid
braided category is clearly defined. For the purpose of this
article, it is more convenient not to emphasize the rigid braided
category in which $(V,c)$ lies. Furthermore, we shall often
consider $(V,c)$ as an object of different categories.
\end{rem}

We refer to \cite{ag} for first properties of Nichols algebras.
We refer to \cite{gfr} for further results, in particular the
freeness result of Nichols algebras over certain Nichols
subalgebras which we state in \ref{fc:ltdf}.

\begin{defn}\label{df:mgr}
Let $\Gm$ be a finite group, $g\in\Gm$,  and
$\rho\in\widehat{\Gm_g}$ an irreducible representation
of $\Gm_g$. We denote by $M(g,\rho)\in\ydskg$ the object
$$M(g,\rho)=\Ind^{\Gm}_{\Gm_g}\rho
	=\k\Gm\otimes_{\k\Gm_g}\rho,\quad
h\ganda(t\otimes x)=ht\otimes x,\quad
\dt(t\otimes x)=tgt^{-1}\otimes (t\otimes x).$$
It has been proved in \cite{dpr,w} and then more generally in
\cite{cibros} in the language of Hopf bimodules (see also
\cite{ag}) that if one takes one element $g$ in each conjugacy
class of $\Gm$ and if $\rho$ runs over all irreducible
non-isomorphic representations of $\Gm_g$, the objects
$M(g,\rho)$ are a full set of irreducible non-isomorphic objects
of $\ydskg$, which turns out to be a semisimple category.
\end{defn}

\begin{defn}\label{df:ddbt}
Let $V=\oplus_{i=1}^{\te}M(g_i,\rho_i)$, where
$\rho_i:\Gm_{g_i}\to\Aut(W_i)$ is an irreducible
representation of $\Gm_{g_i}$.
Let $\{x^i_1,\ldots,x^i_{r_i}\}$ be a basis of $W_i$ and let
$\{h^i_1,\ldots,h^i_{s_i}\}$ be a set of representatives of left cosets
$\Gm/\Gm_{g_i}$. We define 
$t^i_j=h^i_jg_i(h^i_j)^{-1}$ and
$z^i_{jl}=h^i_j\otimes x^i_l\in\Ind^{\Gm}_{\Gm_i} W_i=M(g_i,\rho_i)$.
Then $\{z^i_{jl}\}_{ijl}$ is a basis of $V$ and the braiding is given by
$$c(z^i_{jl}\otimes z^u_{vw})=(t^i_jh^u_v\otimes x^u_w)\otimes z^i_{jl}
=t^i_j\ganda z^u_{vw}\otimes z^i_{jl}.$$
\end{defn}

\begin{defn}
Let $(V,c)$ be a BP, $V=\oplus_iV_i$. We say that this direct
sum is a {\it decomposition of $(V,c)$} (and that $(V,c)$ is
{\it decomposable}) if
$$c|_{V_i\otimes V_j}:V_i\otimes V_j\to V_j\otimes V_i.$$
In this case we denote by $c_{ij}$ the restriction
$c_{ij}=c|_{V_i\otimes V_j}$.
It is straightforward to see that if $(V,c)$ is an RBP and $V=\oplus_iV_i$
is a decomposition of $V$ as a BP, then this decomposition induces a
decomposition of $V^*=\oplus_iV_i^*$ and taking the corresponding
restrictions $c_{ii}$ of $c$ then each $(V_i,c_{ii})$ becomes an RBP.
\end{defn}

If $\Gm$ is a finite group and $V\in\ydskg$ is finite dimensional,
then as noticed in \ref{df:ddbt} $V$ can be decomposed as
$V=\oplus_iM(g_i,\rho_i)$.
However, sometimes it is possible to further decompose some of
the summands, or to give a different decomposition. These decompositions
usually fail to lie in $\ydskg$, but since $(V,c)$ is an RBP any
decomposition gives a direct sum of RBP's.

\begin{defn}\label{df:pttg}
Let $(V,c)$ be a BP. We say that $V$ is {\em of group type}
(or {\em of GT}) if $V$ has a basis $\bcp x1n$ and there exist
$g_i\in\GL(V)\ (i=1,\ldots,n)$ such that
\begin{equation}\label{eq:bat}
c(x_i\otimes x_j)=g_i(x_j)\otimes x_i.
\end{equation}
The terminology, introduced by the referee, comes from the fact
that for $V$ a finite dimensional object of $\ydskg$, one can
find such a basis and maps in $\GL(V)$ (in fact, \ref{df:ddbt}
tells explicitly how to do this when $\Gm$ is finite), and
conversely, a BP which such a property is automatically an object
of $\ydskg$, where $\Gm\subset\GL(V)$ is the group generated
by $\bcp g1n$.
We further say that $(V,c)$ is {\em of finite group type} (or
{\em finite GT}) if the group $\Gm$ generated by $\bcp g1n$
is finite. We say that $(V,c)$ is {\em of abelian group type}
(or {\em abelian GT}) if $\Gm$ is abelian, and we finally say
that $(V,c)$ is {\em of diagonal group type} (or {\em diagonal GT})
if $\Gm$ acts in a diagonal way on $V$. It is clear that
``diagonal" implies ``abelian" and that the latter implies the
former when $\k$ is algebraically closed of characteristic $0$
(as we supposed it to be). It is clear also that ``diagonal"
means that there exist a basis $\bcp x1n$ and a matrix
$(q_{ij})\in(\unid\k)^{n\times n}$ such that
\begin{equation}\label{eq:mybd}
c(x_i\otimes x_j)=q_{ij}x_j\otimes x_i.
\end{equation}
This condition is called ``to come from the abelian case" in
\cite{ag}, since for $V$ in $\ydskg$, $\Gm$ finite abelian,
it is fulfilled. The basis $\bcp x1n$ and the matrix $(q_{ij})$
verifying \ref{eq:mybd} will be called {\em the basis} and
{\em the matrix} of $(V,c)$.

In general, for a BP of group type $(V,c)$, we say that
$\bcp x1n$ is {\em a braided basis} if it verifies \ref{eq:bat}.
The maps $\bcp g1n$ are called {\em the group-likes associated}
to the basis. It can be seen that
\begin{equation}\label{eq:edc}
g_i(x_j)\in V^{g_ig_jg_i^{-1}},
\end{equation}
where we denote $V^g:=\spn\{x_t\ |\ g_t=g\}$.
\end{defn}

It is not hard to prove that a $2$-dimensional
BP of group type is always of abelian GT, and hence if it is
moreover of finite GT, then it is of diagonal GT.

\renewcommand{\theequation}{\arabic{section}.\arabic{equation}}
\section{On the diagonal case}\label{ss:nalg}
We prove an equality in the group algebra of the braid group
$\BB_n$. This equality will give information on the adjoint
action. In particular, the quantum Serre relations will turn
out to be a direct consequence of this equality. We denote by
$e$ the identity element in $\BB_n$ and consider in $\k\BB_n$
the following elements for $1\le i\le j\le n$.
\begin{align*}
U^j_i &= \sgm_j\sgm_{j-1}\cdots\sgm_i, \\
D^j_i &= \sgm_i\sgm_{i+1}\cdots\sgm_j, \\
R^j_i &= (e-\sgm_j\sgm_j)(e-D^j_{j-1}\sgm_j)(e-D^j_{j-2}\sgm_j)
	\cdots(e-D^j_i\sgm_j).
\end{align*}

\begin{lem}\label{lm:qsrg}
We have the following identities
\begin{enumerate}
\item Let $y\in\BB_n\inclu\BB_{n+1}$ and let $x=\iota_{2\_(n+1)}(y)$.
	Then $U^n_1x=yU^n_1$ and \linebreak
	$xD^n_1=D^n_1y$. \label{lm:qsrg1}
\item Let $x\in\BB_{2\_(n+1)}$. Then $xD^n_1U^n_1=D^n_1U^n_1x$. \label{lm:qsrg2}
\item $U^j_1D^n_2\sgm_n=D^n_1\sgm_nU^{j-1}_1$ for $j<n$. \label{lm:qsrg3}
\item $(e+U^1_1+U^2_1+\cdots+U^{n-1}_1-
	D^n_nU^n_1-D^n_{n-1}U^n_1-\cdots-D^n_1U^n_1)R^n_2$\linebreak
	$=R^n_1(e+U^1_1+U^2_1+\cdots+U^{n-1}_1)$ for $n\ge 2$. \label{lm:qsrg4}
\end{enumerate}
\end{lem}

\begin{proof}
\begin{enumerate}
\item Since $\iota_{2\_n}$ is a morphism, we prove it for the generators.
Let $y=\sgm_{j-1}$, $2\le j\le n$. Then $x=\sgm_j$, and we have
\begin{align*}
U^n_1x &=\sgm_n\cdots\sgm_1\sgm_j
	=(\sgm_n\cdots\sgm_{j+1})\sgm_j\sgm_{j-1}\sgm_j
		(\sgm_{j-2}\cdots\sgm_1) \\
&= (\sgm_n\cdots\sgm_{j+1})\sgm_{j-1}\sgm_j\sgm_{j-1}
		(\sgm_{j-2}\cdots\sgm_1)
	=\sgm_{j-1}\sgm_n\cdots\sgm_1 \\
&= yU^n_1,
\end{align*}
and analogously for $D^n_1$.
\item As the previous item, we can prove it for generators $\sgm_j$,
$2\le j\le n$. We have
$$\sgm_jD^n_1U^n_1=D^n_1\sgm_{j-1}U^n_1=D^n_1U^n_1\sgm_j.$$
\item Simply
\begin{align*}
U^j_1D^n_2\sgm_n&=U^j_2D^n_1\sgm_n=\sgm_j\cdots\sgm_2D^n_1\sgm_n=
	\sgm_j\cdots\sgm_3D^n_1\sgm_1\sgm_n=\ldots \\
&\ldots=D^n_1\sgm_{j-1}\cdots\sgm_1\sgm_n
	=D^n_1\sgm_n\sgm_{j-1}\cdots\sgm_1=D^n_1\sgm_nU^{j-1}_1.
\end{align*}
We remark that these equalities are easily understood with drawings.
\item By induction. Notice first that
\begin{align*}
\sgm_i\sgm_{i+1}\sgm_{i+1}\sgm_i\sgm_{i+1}\sgm_{i+1}
	&=\sgm_i\sgm_{i+1}\sgm_i\sgm_{i+1}\sgm_i\sgm_{i+1}
	=\sgm_{i+1}\sgm_i\sgm_{i+1}\sgm_i\sgm_{i+1}\sgm_i \\
	&=\sgm_{i+1}\sgm_{i+1}\sgm_i\sgm_{i+1}\sgm_{i+1}\sgm_i.
\end{align*}
For $n=2$ we have
\begin{align*}
(e&+\sgm_1-\sgm_2\sgm_2\sgm_1-\sgm_1\sgm_2\sgm_2\sgm_1)(e-\sgm_2\sgm_2) \\
&= e-\sgm_2\sgm_2+\sgm_1-\sgm_1\sgm_2\sgm_2-\sgm_2\sgm_2\sgm_1+
	\sgm_2\sgm_2\sgm_1\sgm_2\sgm_2-\sgm_1\sgm_2\sgm_2\sgm_1
	+\sgm_1\sgm_2\sgm_2\sgm_1\sgm_2\sgm_2 \\
&= e-\sgm_2\sgm_2+\sgm_1-\sgm_1\sgm_2\sgm_2-\sgm_2\sgm_2\sgm_1+
	\sgm_2\sgm_2\sgm_1\sgm_2\sgm_2-\sgm_1\sgm_2\sgm_2\sgm_1
	+\sgm_2\sgm_2\sgm_1\sgm_2\sgm_2\sgm_1 \\
&= (e-\sgm_2\sgm_2-\sgm_1\sgm_2\sgm_2+\sgm_2\sgm_2\sgm_1\sgm_2\sgm_2)
	(e+\sgm_1) \\
&= (e-\sgm_2\sgm_2)(e-\sgm_1\sgm_2\sgm_2)(e+\sgm_1)
\end{align*}
Suppose the equality is true for $n$. Then
\begin{align*}
(e&+U^1_1+U^2_1+\cdots+U^n_1-D^{n+1}_{n+1}U^{n+1}_1
	-D^{n+1}_nU^{n+1}_1-\cdots-D^{n+1}_1U^{n+1}_1)R^{n+1}_2 \\
&= R^{n+1}_2+(U^1_1+U^2_1+\cdots+U^n_1-D^{n+1}_{n+1}U^{n+1}_1-\cdots
	-D^{n+1}_2U^{n+1}_1)R^{n+1}_2 \\
&\qquad -D^{n+1}_1U^{n+1}_1R^{n+1}_2 \\
&= R^{n+1}_2+(e+U^2_2+\cdots+U^n_2-D^{n+1}_{n+1}U^{n+1}_2-\cdots
	-D^{n+1}_2U^{n+1}_2)\sgm_1R^{n+1}_2 \\
&\qquad -D^{n+1}_1U^{n+1}_1R^{n+1}_2 \\
&= R^{n+1}_2+(e+U^2_2+\cdots+U^n_2-D^{n+1}_{n+1}U^{n+1}_2-\cdots
	-D^{n+1}_2U^{n+1}_2)\sgm_1R^{n+1}_3 \\
&\qquad \times(e-D^{n+1}_2\sgm_{n+1})-R^{n+1}_2D^{n+1}_1U^{n+1}_1
		\mbox{ (by \ref{lm:qsrg2})} \\
&= R^{n+1}_2+(e+U^2_2+\cdots+U^n_2-D^{n+1}_{n+1}U^{n+1}_2-\cdots
	-D^{n+1}_2U^{n+1}_2)R^{n+1}_3\sgm_1 \\
&\qquad \times(e-D^{n+1}_2\sgm_{n+1})-R^{n+1}_2D^{n+1}_1U^{n+1}_1
		\mbox{ ($R^{n+1}_3\in\k\BB_{3\_n+2}$)} \\
&= R^{n+1}_2+\iota_{2\_(n+2)}((e+U^1_1+\cdots+U^{n-1}_1-D^n_nU^n_1-\cdots
	-D^n_1U^n_1)R^n_2)\sgm_1 \\
&\qquad \times (e-D^{n+1}_2\sgm_{n+1})-R^{n+1}_2D^{n+1}_1U^{n+1}_1 \\
&= R^{n+1}_2+\iota_{2\_(n+2)}(R^n_1(e+U^1_1+U^2_1+\cdots+U^{n-1}_1))
	\sgm_1(e-D^{n+1}_2\sgm_{n+1}) \\
&\qquad -R^{n+1}_2D^{n+1}_1U^{n+1}_1 \mbox{ (by IH)} \\
&= R^{n+1}_2+R^{n+1}_2(e+U^2_2+U^3_2+\cdots+U^n_2)
	\sgm_1(e-D^{n+1}_2\sgm_{n+1})
	-R^{n+1}_2D^{n+1}_1U^{n+1}_1 \\
&= R^{n+1}_2[e+\sgm_1+U^2_2\sgm_1+\cdots+U^n_2\sgm_1-
	\sgm_1D^{n+1}_2\sgm_{n+1} \\
&\qquad -U^2_2\sgm_1D^{n+1}_2\sgm_{n+1}
	-\cdots-U^n_2\sgm_1D^{n+1}_2\sgm_{n+1}-D^{n+1}_1U^{n+1}_1] \\
&= R^{n+1}_2(e+U^1_1+U^2_1+\cdots+U^n_1-D^{n+1}_1\sgm_{n+1}
	-U^2_1D^{n+1}_2\sgm_{n+1}-\cdots \\
&\qquad -U^n_1D^{n+1}_2\sgm_{n+1}-D^{n+1}_1\sgm_{n+1}U^n_1) \\
&= R^{n+1}_2(e+U^1_1+\cdots+U^n_1-D^{n+1}_1\sgm_{n+1}
	-D^{n+1}_1\sgm_{n+1}U^1_1-\cdots-D^{n+1}_1\sgm_{n+1}U^{n-1}_1 \\
&\qquad -D^{n+1}_1\sgm_{n+1}U^n_1) \mbox{ (by \ref{lm:qsrg3})} \\
&= R^{n+1}_2(e-D^{n+1}_1\sgm_{n+1})(e+U^1_1+U^2_1+\cdots+U^n_1) \\
&=R^{n+1}_1(e+U^1_1+U^2_1+\cdots+U^n_1).
\end{align*}
\end{enumerate}
\end{proof}

\begin{prop}\label{pr:qsrg}
$\SSS^{n+1}(e-U^n_1)(e-U^n_2)\cdots(e-U^n_n)=R^n_1\SSS^{n|1}$
	for $n\ge 1$.
\end{prop}
\begin{proof}
By induction, for $n=1$ being
$$\SSS^2(e-\sgm_1)=(e+\sgm_1)(e-\sgm_1)=e-\sgm_1\sgm_1.$$
If the equality holds for $n$, we have
\begin{align*}
\SSS^{n+2} & (e-U^{n+1}_1)(e-U^{n+1}_2)\cdots(e-U^{n+1}_{n+1}) \\
&= \SSS^{n+2}(e-U^{n+1}_2)\cdots(e-U^{n+1}_{n+1})
	-\SSS^{n+2}U^{n+1}_1(e-U^{n+1}_2)\cdots(e-U^{n+1}_{n+1}) \\
&= \TTT_{1,n+1}\SSS^{1|n+1}(e-U^{n+1}_2)\cdots(e-U^{n+1}_{n+1}) \\
&\qquad -\TTT_{n+1,1}\SSS^{n+1|1}
	U^{n+1}_1(e-U^{n+1}_2)\cdots(e-U^{n+1}_{n+1}) \\
&= \TTT_{1,n+1}\SSS^{1|n+1}(e-U^{n+1}_2)\cdots(e-U^{n+1}_{n+1}) \\
&\qquad -\TTT_{n+1,1}U^{n+1}_1\SSS^{1|n+1}
	(e-U^{n+1}_2)\cdots(e-U^{n+1}_{n+1})
		\mbox{ (by \ref{lm:qsrg},\ref{lm:qsrg1})} \\
&= (\TTT_{1,n+1}-\TTT_{n+1,1}U^{n+1}_1)\SSS^{1|n+1}
	(e-U^{n+1}_2)\cdots(e-U^{n+1}_{n+1}) \\
&= (\TTT_{1,n+1}-\TTT_{n+1,1}U^{n+1}_1)\iota_{2\_n+2}(\SSS^{n+1}
	(e-U^n_1)\cdots(e-U^n_n)) \\
&= (\TTT_{1,n+1}-\TTT_{n+1,1}U^{n+1}_1)\iota_{2\_n+2}
	(R^n_1\SSS^{n|1}) \mbox{ (by IH)} \\
&= (\TTT_{1,n+1}-\TTT_{n+1,1}U^{n+1}_1)(R^{n+1}_2\SSS^{1|n|1}) \\
&= (e+U^1_1+U^2_1+\cdots+U^{n+1}_1-U^{n+1}_1-D^{n+1}_{n+1}U^{n+1}_1
	-\cdots-D^{n+1}_1U^{n+1}_1) \\
&\qquad \times R^{n+1}_2\SSS^{1|n|1} \\
&= (e+U^1_1+\cdots+U^n_1-D^{n+1}_{n+1}U^{n+1}_1
	-\cdots-D^{n+1}_1U^{n+1}_1)R^{n+1}_2\SSS^{1|n|1} \\
&= R^{n+1}_1(e+U^1_1+\cdots+U^n_1)
	\SSS^{1|n|1} \mbox{ (by \ref{lm:qsrg},\ref{lm:qsrg4})} \\
&= R^{n+1}_1(\TTT_{1,n}|e_1)(\SSS^{1|n}|e_1)
	\quad\mbox{($e_1$ the unique element of $\BB_1$)} \\
&= R^{n+1}_1(\SSS^{n+1}|e_1) \\
&= R^{n+1}_1\SSS^{(n+1)|1}.
\end{align*}
\end{proof}

\begin{defn}
Let $X$ be a vector space, $0\neq x\in X$ and $T\in\End(X)$
such that $T^n(x)=0$ for some $n$. The {\em nilpotency order}
of $T$ in $x$ is the number $r+1$ such that $T^r(x)\neq 0$
and $T^{r+1}(x)=0$.
\end{defn}

\begin{defn}
Let $R$ be a braided Hopf algebra and $x\in R$. We denote
by $\Ad_x\in\End(R)$ the (braided) adjoint of $x$, i.e.,
$$\Ad_x(y)=m(\com x1\otimes m(c(\Ss(\com x2)\otimes y)),$$
where $c$ is the braiding of $R$ and $m$ the multiplication
map. Thus, if $x$ is a primitive element, we have
$$\Ad_x(y)=xy-m(c(x\otimes y)).$$
\end{defn}

The following result can be found in \cite{ros}, but, as said in
the introduction, we can obtain it directly from \ref{pr:qsrg}.
\begin{thm}\label{th:qsrg}
Let $x_1,x_2\in V$ be such that for some matrix $(q_{ij})$
the braiding verifies
$c(x_i\otimes x_j)=q_{ij}x_j\otimes x_i$ for $i,j=1,2$.
Let $N=N(q_{11})$ (see the beginning of \ref{ss:bp}) and
let $t$ be the least non-negative integer (if it exists) such
that $q_{11}^tq_{12}q_{21}=1$. Let $r=\min\{t,N-1\}$. Then
the nilpotency order of $\Ad_{x_1}$ on $x_2$ is exactly $r+1$.
\end{thm}
\begin{proof}
Let $n\in\NN$. We consider $\Ad_{x_1}^n(x_2)=\SSS^{n+1}
(\widetilde{\Ad}_{x_1}^n(x_2))\in\toba V$, where
$\widetilde{\Ad}$ is the adjoint in $AV$. It is easy to
see by induction that this element is the action of the
left hand side of \ref{pr:qsrg} on $x_1^nx_2$. Now, the
right hand side acts on $x_1^nx_2$ as a scalar, namely
$$(1-q_{12}q_{21})(1-q_{11}q_{12}q_{21})
(1-q_{11}^2q_{12}q_{21})\cdots(1-q_{11}^{n-1}q_{12}q_{21})
	\fsq n{q_{11}}.$$
The assertion follows at once.
\end{proof}

\begin{defn}\label{df:der}
Let $(V,c)$ be an RBP. For $y\in V^*$, we denote by $\ptl_y$ the operator
$$\ptl_y=(\id\otimes y)\circ\Dt^{i-1,1}:\tobag iV\to\tobag {i-1}V.$$
If $B=\{x_1,\ldots,x_n\}$ is a basis of $V$, let $\{\ptl_1,\ldots,\ptl_n\}$
be the basis of $V^*$ dual to $B$.
We denote also by $\ptl_j,\ j=1,\ldots,n$ the operator $\ptl_{\ptl_j}$.
We warn that for the case $V\in\yd$, these are usually not morphisms
in $\yd$, but only linear maps.
\end{defn}

We state some facts which we shall use. They are consequences
of results in \cite{gfr}.
\begin{thm}\label{fc:ltdf}
Let $V=\oplus_{i=1}^nM(g_i,\rho_i)=\oplus_iV_i$, where
$\rho_i\in\widehat{\Gm_{g_i}}$. Let $q_i=\rho_i(g_i)$ (since $g_i$
lies in the center of $\Gm_{g_i}$ and $\rho_i$ is irreducible,
$\rho_i(g_i)$ is a scalar) and let $N_i=N(q_i)$.
\begin{enumerate}
\item\label{fc:ltdf1}
Take a basis $\{z^i_{jl}\}$ of $V$ as in \ref{df:ddbt}
and $\{\ptl^i_{jl}\}$ its dual basis. Let $i,j,l$ be
fixed, $K=\ker(\ptl^i_{jl})\subset\toba V$, and let
$W=\spn(x^i_{jl})$. Then for those $i$ such that $N_i<\infty$
the multiplication map
$K\otimes\toba W\to\toba V$ is an isomorphism. In
particular, if $\{y_1,\ldots,y_k\}$ is a basis of $K$,
then
$$\{y_r(z^i_{jl})^s\ |\ 1\le r\le k,\ 0\le s<N_i\}$$
is a basis of $\toba V$.
\item\label{fc:ltdf2}
If $V$ is of finite Nichols rank, then
$$N_i\ |\ \dim\toba{V_i}\ |\ \prod_i\dim\toba{V_i}\ |\ \dim\toba V.$$
Furthermore, $\prod_i\dim\toba{V_i}=\dim\toba V$ iff
$c_{ij}=\id_{V_i\otimes V_j}\ \forall i\neq j$, where
$c_{ij}=c|_{V_i\otimes V_j}:V_i\otimes V_j\to V_j\otimes V_i$.
\item\label{fc:ltdf3}
$N_i^{\dim V_i}\le\dim\toba{V_i}\ \forall i$. Furthermore, if
$\bcp{z^i}1{d_i}$ is a basis of $V_i$, then the set
$$\{(z^1_1)^{i^1_1}\cdots(z^1_{d_1})^{i^1_{d_1}}(z^2_1)^{i^2_1}
\cdots(z^n_{d_n})^{i^n_{d_n}}\ |\ 0\le i^j_k<N_j\}$$
is linearly independent in $\toba V$.
\end{enumerate}
\end{thm}

This allows to give a relation between the dimension and the Nichols
rank of a BP of finite group type:
\begin{defn}\label{df:rd}
Let $n=\prod_{i=1}^Rp_i^{v_i}$, where $p_i$ are distinct primes, $v_i>0$
and $p_i>p_1$ for $i>1$. We define $r(n)=\log_{p_1}(n)$.
\end{defn}

\begin{rem}\label{rm:csm}
Applying \ref{fc:ltdf} part \ref{fc:ltdf1} together with part
\ref{fc:ltdf3}, we get that if \linebreak
$\dim\toba V=n$, then $d=\dim V\le r(n)$.
Applying \ref{fc:ltdf} part \ref{fc:ltdf2}, we get that if
$V=\oplus_{i=1}^{\te}V_i=\oplus_iM(g_i,\rho_i)\in\ydskg$
and $n=\dim(\toba V)$, $n=\prod_{i=1}^Rp_i^{v_i}$, then
$\te\le\sum_iv_i$.
\end{rem}

\begin{defn}\label{df:qls}
If $(V,c)$ is of diagonal GT with basis $\bcp x1n$ and matrix
$(q_{ij})$, we say, following \cite{as2}, that $\toba V$ is a
{\em Quantum Linear Space} (or QLS) if $q_{ij}q_{ji}=1$ for
$i\neq j$. It is proved in \cite{as2} (and it is also a consequence
of \ref{fc:ltdf} part \ref{fc:ltdf2}), that in this case
$\dim\toba V=\prod N(q_{ii})$.
\end{defn}

In \cite{as3} the authors associate to some Nichols algebras
of diagonal group type a Frobenius--Lusztig kernel with
the same dimension. This association can be done whenever
$q_{ij}q_{ji}=q_{ii}^{a_{ij}}\ \forall i,j$ for a matrix $(a_{ij})$
with entries in $\ZZ$, where $(q_{ij})$ is the matrix of the
braiding. The following proposition deals with general rank-$2$
Nichols algebras.

\begin{prop}\label{pr:mdr2}
Let $(V,c)$ be a $2$ dimensional BP of finite GT (and thus of
diagonal GT), with basis $\{x_1,x_2\}$ and matrix $(q_{ij})$,
and let $N_j=N(q_{jj})$ for $j=1,2$.
As in \ref{th:qsrg}, let $r+1$ be nilpotency order of
$\Ad_{x_2}$ on $x_1$. For $1\le i\le r$ let
$M_i=N(q_{11}(q_{12}q_{21})^iq_{22}^{i^2})$. Then
$$\dim\toba V\ge N_1N_2\prod_{1\le i\le r}M_i.$$
Furthermore, suppose that the nilpotency order of
$\Ad_{x_1}$ on $x_2$ is $2$. We have:
\begin{enumerate}
\item If $r=1$ then the equality holds.
\item If $r=2$, and $N(q_{11})\neq 2$ or $N(q_{22})\neq 3$ then
the equality holds.
\item If $r=2$, $N(q_{11})=2$ and $N(q_{22})=3$ then the equality
holds if and only if
$$q_{12}q_{21}=-1,\mbox{ or }q_{12}q_{21}=q_{22},
\mbox{ or }q_{12}q_{21}=-q_{22}.$$
\end{enumerate}
\end{prop}
\begin{proof}
Let $z_1=\Ad_{x_2}(x_1)=x_2x_1-q_{21}x_1x_2$,
and let us define inductively \linebreak
$z_{i+1}=\Ad_{x_2}(z_i)=x_2z_i-q_{21}q_{22}^iz_ix_2$.
Taking on $\toba V$ the bi-degree given by
$\deg(x_1)=(1,0)$, $\deg(x_2)=(0,1)$, we have
$z_i\in\tobag {1,i}V$. Furthermore, it is proved in
\ref{th:qsrg} that $z_r\neq 0$ and $z_{r+1}=0$. It is
immediate to see by induction (it follows also from
\ref{pr:qsrg}) that $\ptl_{x_2}(z_i)=0$, and
$\ptl_{x_1}(z_i)=d_{12}^1d_{12}^2\cdots d_{12}^ix_2^i$,
where $d_{12}^j=(1-q_{12}q_{21}q_{22}^{j-1})$. Thus,
$\ptl_{x_2}^i\ptl_{x_1}(z_i)\neq 0$ if $i\le r$.
We shall prove that the set
$\{x_1^{n_1}z_1^{m_1}\cdots z_r^{m_r}\ |
	\ 0\le n_1<N_1,\ 0\le m_i<M_i\}$
is linearly independent. To see this, we compute
\begin{align*}
\ptl_{x_2}^i\ptl_{x_1}(z_i^m) &=
	\sum_{j=0}^{m-1}q_{12}^{ji}q_{11}^jd_{12}^1\cdots d_{12}^j
		\ptl_{x_2}^i(z_i^{m-1-j}x_2^iz_i^j) \\
	&= d_{12}^1\cdots d_{12}^j\sum_{j=0}^{m-1}
	q_{12}^{ji}q_{11}^jq_{21}^{ji}q_{22}^{ji^2}
	\fsq i{q_{22}}z_i^{m-1} \\
	&=d_{12}^1\cdots d_{12}^j\nsq m{q_{11}q_{12}^iq_{21}^iq_{22}^{i^2}}
	z_i^{m-1},
\end{align*}
from where $z_i^m\neq 0$ if
$m<N(q_{11}q_{12}^iq_{21}^iq_{22}^{i^2})=M_i$.
It is clear also that if $j>i$ then
$\ptl_{x_2}^j\ptl_{x_1}(z_i^m)=0$. Suppose
inductively that for $s<r$, the set
$$\{x_1^{n_1}z_1^{m_1}\cdots z_s^{m_s}\ |
	\ 0\le n_1<N_1,\ 0\le m_i<M_i\}$$
is linearly independent. Then we have that if there
were a linear combination
$$\sum\al_{n_1,m_1,\ldots,m_{s+1}}
x_1^{n_1}z_1^{m_1}\cdots z_{s+1}^{m_{s+1}}=0,$$
we would get, applying
$(\ptl_{x_2}^{s+1}\ptl_{x_1})^{M_{s+1}-1}$, a linear combination
$$\sum\tilde\al_{n_1,m_1,\ldots,m_s,M_{s+1}-1}
	x_1^{n_1}z_1^{m_1}\cdots z_s^{m_s}=0,$$
where the $\tilde\al$'s are the $\al$'s multiplied
by non-zero factors. By the inductive  \linebreak
hypothesis, we have that $\al_{n_1,m_1,\ldots,m_s,M_{s+1}-1}=0$.
We apply now $[\ptl_{x_2}^{s+1}\ptl_{x_1}]^{M_{s+1}-2}$
and get that $\al_{n_1,m_1,\ldots,m_s,M_{s+1}-2}=0$.
Continuing in this way, we get that \linebreak
$\al_{n_1,m_1,\ldots,m_s,m_{s+1}}=0$ for all $m_{s+1}$,
proving the inductive thesis. Now the 
assertion follows at once noting that this set lies in
$\ker(\ptl_{x_2})$, and then by \ref{fc:ltdf} part
\ref{fc:ltdf1} the set
\begin{equation}\label{eq:sbpbw}
\{x_1^{n_1}z_1^{m_1}\cdots z_r^{m_r}x_2^{n_2}\ |
	\ 0\le n_i<N_i,\ 0\le m_i<M_i\}
\end{equation}
is linearly independent in $\toba V$.

Now, for the nilpotency order of $\Ad_{x_1}$ on $x_2$ to be $2$
it can happen that $q_{11}=-1$ or that $q_{11}q_{12}q_{21}=1$.
In the first case we have $x_1^2=0$, and then
$$z_1x_1=-q_{21}x_1x_2x_1=q_{11}q_{21}z_1x_1.$$
In the second case it is easy to see using derivations that
we also have
$$z_1x_1=q_{11}q_{21}z_1x_1,$$
for $\ptl_{x_2}$ annihilates and $\ptl_{x_1}$ gives
the same element in $\toba V$ when applied to both sides of
the equality.

We consider now the three cases in the statement.
\begin{enumerate}
\item 
We have $0=\Ad_{x_2}(z_1)=x_2z_1-q_{21}q_{22}z_1x_2$,
whence $x_2z_1=q_{21}q_{22}z_1x_2$. These equalities
show how to write any monomial in $\toba V$ as a
combination of the elements of the set
$\{x_1^{n_1}z_1^{m_1}x_2^{n_2}\}$, from where
$\dim\toba V\le N_1N_2M_1$ and the equality follows.

\item
We have $z_2=\Ad_{x_2}(z_1)=x_2z_1-q_{21}q_{22}z_1x_2$,
and hence
$$0=\Ad_{x_2}(z_2)=x_2z_2-q_{21}q_{22}^2z_2x_2.$$
It can be seen, using derivations, that
$z_2x_1=q_{21}^2q_{11}x_1z_2+(q_{11}-q_{22})q_{21}z_1^2$.
Furthermore, it can also be seen that
$z_2z_1=-q_{21}z_1z_2$. We conclude now as in the previous
case, since these equalities show how to write any monomial
as a combination of those of \eqref{eq:sbpbw}.

\item The same computations as in the previous case hold
here, except for the relation between $z_2z_1$ and
$z_1z_2$. We are looking for the condition on the
matrix for \eqref{eq:sbpbw} to be a basis of $\toba V$.
For this to happen, we must be able to write $z_2z_1$ as a linear
combination of those elements; but now, the bi-degree
of $z_2z_1$ is $(3,2)$ and, moreover, it lies in the
kernel of $\ptl_{x_2}$. By \ref{fc:ltdf} part \ref{fc:ltdf1},
we have to write $z_2z_1$ as a combination of the
elements in \eqref{eq:sbpbw} of bi-degree $(3,2)$ which
also lie in $\ker\ptl_{x_2}$, but there is only one such
element: $z_1z_2$. Hence we should have $z_2z_1=\lmb z_1z_2$.
Taking $\ptl_{x_1}$ in both sides, we are led to the equations
\begin{align*}
&-\lmb(1-t)tq_{12}q_{22}^2=(1-t)(1+t-q_{22}t^2), \\
&\lmb(1-t)(1-qt)=-(1-t)(1-qt)tq_{21}q_{22}^2,
\end{align*}
for $t=q_{12}q_{21}$. Since $(1-t)\neq 0$ and $(1-qt)\neq 0$,
these equations are equivalent to
$$0=1+t-q_{22}t^2-q_{22}t^3=q_{22}(1+t)(q_{22}-t)(q_{22}+t),$$
from where the conditions we stated are clear.
\end{enumerate}
\end{proof}

\begin{rem}\label{rm:cpnc}
For $(V,c)$ a BP of diagonal GT, let the matrix $a_{ij}$
be defined by $a_{ii}=2$ and $a_{ij}=1-d_{ij}$ for
$i\neq j$, with $d_{ij}$ the nilpotency order of $\Ad_{x_i}$
on $x_j$. This is a generalized Cartan matrix. The proposition
shows that if $(a_{ij})$ is of type $A_2$ then a PBW basis can
be constructed for the Nichols algebra $\toba V$
in the same way as for Frobenius--Lusztig kernels of type $A_2$,
though these matrices are not of Cartan type in general.
It is natural to ask to what extent this method allows
to compute the Nichols algebras as generalizations
of Frobenius--Lusztig kernels. In part 3 we show that this
does not give the general answer, and we find the conditions
for this method to work in the case of type $B_2$.

Concrete examples of these algebras over specific groups
can be found in \linebreak
\cite[\S 3]{n}. There, the author presents
algebras over abelian groups, all of which can be computed
with \ref{pr:mdr2}. For instance, if $\Gm=C_4$ the
cyclic group of order $4$ with generator $\sgm$, take
$V=M(\sgm,\chi)\oplus M(\sgm^2,\chi)$, where $i=\chi(\sgm)$
is a primitive fourth root of unity. Then the matrix
of $c$ is $\mdpd {-1}i{-1}i$. Here $N_1=2$, $N_2=4$,
$r=1$ and the nilpotency order of $\Ad_{x_1}$ on $x_2$ is $2$.
Then $\toba V$ is of type $A_2$, and $M_1=N(-1)=2$, whence
$\dim\toba V=16$.  The pointed Hopf algebra $\toba V\#\k C_4$
is hence $64$-dimensional.
As another example, borrowed from \cite{n},
take $\Gm=C_6$ with generator $\sgm$, and
take $V=M(\sgm,\chi^3)\oplus(\sgm,\chi^2)$, where
$N(\chi(\sgm))=6$, i.e. $\chi$ is a generator of $\widehat\Gm$.
Then the matrix of $c$ is $\mdpd {-1}\omega{-1}\omega$, where
$\omega=\chi^2(\sgm)$ is a third root of unity. This braiding
is of type $B_2$ and fits in the part 3 of the proposition,
with $q_{12}q_{21}=-q_{22}$. We have $N_1=2$, $N_2=3$, $M_1=3$,
$M_2=2$, and hence $\toba V$ is $36$ dimensional.
More examples can be found in \cite{g32}.
\end{rem}

\section{On the non-diagonal case}
We give some necessary conditions a Yetter--Drinfeld
module over a finite group algebra must verify to produce
a finite dimensional Nichols algebra. As said in \ref{df:pttg},
such a module is of group type. We want, however, to keep
track of the structure of the group.

\begin{lem}\label{lm:csgr}
Let $\Gm$ be a finite group, $g\in\Gm$ and
$\rho\in\widehat{\Gm_g}$ an irreducible representation
of $\Gm_g$. Let $q=\rho(g)$ (as noticed in \ref{fc:ltdf},
it is a scalar) and let $N=N(q)$. Let $V=M(g,\rho)$
and suppose $V$ is of finite Nichols rank. Then
\begin{enumerate}
\item if $\deg\rho\ge 3$ then $N=2$ (i.e. $q=-1$);
\item if $\deg\rho=2$ then $2\le N\le 3$
(i.e. $q=-1$ or it is a third root of unity).
\end{enumerate}
\end{lem}
\begin{proof}
Let $W$ be the space of the representation $\rho$.
For $x\in W$ we denote also by $x$ the element
$1\otimes x\in\Ind^{\Gm}_{\Gm_g}W=V$. Let $\bcp x1d$
be a basis of $W$.
We have $\dt(x)=g\otimes x\ \forall x\in W$, and hence
$c(x_i\otimes x_j)=g\ganda x_j\otimes x_i=qx_j\otimes x_i$.
We have an inclusion $\toba W\inclu\toba V$ (in fact,
$\toba V$ is free as a right $\toba W$-module by \linebreak
\cite[Th. 3.11]{gfr}) and hence $\toba W$ is finite
dimensional. Now, $W$ is of diagonal GT with matrix
$(q_{ij})$, $q_{ij}=q\ \forall i,j$, and thus it is
of Cartan type. Hence, the Cartan matrix of
this braiding has $2$ in the main diagonal and $2-N$
outside the main diagonal. Since by a result of Lusztig
(see \cite[Thm. 3.1]{as-cm}) this matrix must correspond
to a finite datum, the condition on $N$ follows.
\end{proof}

As a consequence of this, we see that Nichols algebras which
do not come from the abelian case are usually infinite
dimensional. In the same vein, we have the following result
which appeared as a conjecture in a former version of \cite{as-cm}.
\begin{prop}
If $\Gm$ has odd order, then there are finitely many modules
in $\ydskg$ of finite Nichols rank.
\end{prop}
\begin{proof}
Since $\ydskg$ has a finite number of simple modules up to
isomorphism, it is enough to see that each simple module
$M(g,\rho)$ can appear no more than $n$ times ($n$ depending
on $M(g,\rho)$) in a module $V$ for it to be of finite Nichols
rank. Let then $W$ be the space affording $\rho$, $0\neq x\in W$
and $q$ such that $g\ganda x=qx$. Suppose $V=n\cdot M(g,\rho)+\cdots$
($n>0$). Then we take $Y=\spn\bcp x1n$, where $x_i$ is the element
of the i-th copy of $M(g,\rho)$ corresponding to $x$. It is
immediate to see that $Y$ is of diagonal GT with matrix
$(q_{ij})$, $q_{ij}=q\ \forall i,j$. Now, $q\neq 1$ since
$M(g,\rho)$ is of finite Nichols rank (since so is $V$),
and $q\neq -1$ since $\Gm$ is of odd order. The same argument
as in the previous lemma tells that $1\le n\le 2$ if $N(q)=3$
and $n=1$ if $N(q)>3$.
\end{proof}
Without the odd-order assumption the result is known to be
false (for instance, for $\Gm=C_2$).

\medskip
Let now $\Gm$ be a finite group and $\orb_g$ the conjugacy
class of $g\in\Gm$. Let \linebreak
$\rho:\Gm_g\to\Aut(W)$ be an irreducible representation, and take
a basis of $M(g,\rho)$ as follows: let $\{h_1,\ldots,h_s\}$ be a set
of representatives of left cosets of $\Gm/\Gm_g$ and let
$\{x_1,\ldots,x_r\}$ be a basis of $W$. We consider the basis
$\{z_{jl}=h_j\otimes x_l\}$. By the definition of the braiding, it is
useful to take into account the way in which $\orb_g$ acts on itself
by conjugation. To do this, we name $t_i=h_igh_i^{-1}$ and define the
morphism
$$f_g:\Gm\to\SS_s,\quad f_g(k)(i)=j\ \mbox{if}\ kt_ik^{-1}=t_j.$$
Notice that $f_g(k)$ fixes $i$ iff $k$ commutes with $t_i$; in particular
$f_g(t_i)(i)=i$ for all $i$.

\begin{rem} \label{rm:osf}
Let $f:G\to H$ be a morphism between two groups. For
$g\in G$, the conjugacy class $\orb_g$ has cardinality
the index $s=[G:G_g]$. Let $\tilde G=f^{-1}H_{f(g)}$.
It is straightforward to see that the set $f(\orb_g)$
has cardinality the index $n=[G:\tilde G]$ and that for
any $h\in f(\orb_g)$ the fiber $f^{-1}(h)\cap\orb_g$ has
cardinality the index $m=[\tilde G:G_g]$, whence the
cardinality of the orbit $\orb_g$ can be factorized by $s=nm$.
\end{rem}

\begin{rem}\label{rm:saq}
Let $V=M(g,\rho)$ where $\rho:\Gm_g\to\Aut(W)$ and let
$x\in W$. Let $q=\rho(g)$ (as noticed in \ref{fc:ltdf},
$g$ acts by a scalar, on $W$). In particular, \linebreak
$c(x\otimes x)=g\ganda x\otimes x=qx\otimes x$.
Let $t\notin\Gm_g$ and $y=t\ganda x$. Thus
$\dt(y)=tgt^{-1}\otimes y$ and $tgt^{-1}$
acts on $y$ by $q$, since
$$tgt^{-1}\ganda y=tgt^{-1}t\ganda x
	=tg\ganda x=qt\ganda x=qy.$$
Thus, $c(y\otimes y)=qy\otimes y$.
\end{rem}

The preceding paragraphs suggest a different approach for
classifying braided pairs of finite group type, as follows:
\begin{defn}\label{df:chcc}
Let $X$ be a finite set and $\trid:X\times X\to X$ a function.
We say that $(X,\trid)$ is a {\em crossed set} if
\begin{enumerate}
\item for each $i\in X$, the function $i\trid\bullet:X\to X$,
$j\mapsto i\trid j$, is a bijection,\label{ccc1}
\item $i\trid i=i\ \forall i\in X$,\label{ccc2}
\item $j\trid i=i$ whenever $i\trid j=j$, and\label{ccc3}
\item $i\trid(j\trid k)=(i\trid j)\trid(i\trid k)\ \forall i,j,k\in X$.\label{ccc4}
\end{enumerate}
The example we come from is $X$ being a union of conjugacy classes
of a finite group $\Gm$ and $\trid$ the conjugation, $i\trid j=iji^{-1}$.

A crossed set provides a set-theoretical solution to the Braid Equation by
\begin{equation}\label{eq:calt}
c:X\times X\to X\times X,\quad c(i,j)=(i\trid j,i).
\end{equation}
We remark that $c$ being a solution to the Braid Equation is equivalent
to the condition \ref{ccc4} in the above definition. The condition
\ref{ccc1} is given to guarantee the bijectivity of $c$, and the
conditions \ref{ccc2} and \ref{ccc3} are necessary conditions for $X$
to be ``injective", in the language of \cite{s}. With respect to our
intentions of producing braided pairs of finite GT, the conditions
\ref{ccc2} and \ref{ccc3} are harmless:
it is not hard to see that if $(X,\trid)$ verifies \ref{ccc1} and
\ref{ccc4} and $f\in C^2(X)$ is a $2$-cocycle then there exist
$X',\ \trid',\ f'$ such that $(X',\trid')$ is a crossed set and
$(X',c^{f'})$ is isomorphic to $(X,c^f)$ (see below for the notation).

Notice that any set-theoretical solution to the Braid Equation
$c$ gives rise to a (set-theoretic) representation of the braid
group $\BB_n$, namely, $X^n=X\times\cdots\times X$ ($n$ times),
where $\sgm_i$ acts by $c$ in the coordinates $i$ and $i+1$.

Although \cite{s,lyz} consider more general set-theoretical solutions,
in both articles the authors prove that any bijective solution
$(i,j)\mapsto ({}^ij,i^j)$ is, in some sense, equivalent to a
solution as in \eqref{eq:calt}. Specifically, for any solution
$\tilde c$ there exists a solution $c$ as in \eqref{eq:calt}
and a collection of bijective maps $J_n:X^n\to X^n$ which intertwine
the representations of $\BB_n$ associated to $\tilde c$ and $c$.

Let $(X,\trid)$ be a crossed set and $n\ge 0$. We define
$C^n(X):=\{f:X^n\to\unid\k\}$ ($X^0$ is a singleton).
Let $\dt^n:C^n(X)\to C^{n+1}(X)$ be given by
\begin{align*}
\dt^0&=0 \\
\dt^n(f)(x_0,\ldots,x_n)&=\prod_{i=0}^{n-1}
	f(x_0,\ldots,x_{i-1},x_{i+1},\ldots,x_n)^{(-1)^i} \\
&\hspace{2cm} \times f(x_0,\ldots,x_{i-1},x_i\trid x_{i+1},\ldots,x_i\trid x_n)^{(-1)^{i+1}}.
\end{align*}
It is easy to prove that $(C^\bullet(X),\dt)$ is a cochain complex (we remark
that an analog definition for a general set-theoretical solution of the
Braid Equation does not yield a complex). We define then the cohomology
$H^n(X)=H^n(C^\bullet(X),\dt)=\ker(\dt^n)/\im(\dt^{n-1})$.
\end{defn}

Let $f\in C^2(X)$. We define a map $c^f:\k X\otimes\k X\to\k X\otimes\k X$ by
$$c^f(i\otimes j)=f(i,j)i\trid j\otimes i,$$
and extend it linearly. We denote by $\bdng(X)$ the subset
$$\bdng(X)=\{f\in C^2(X)\ |\ c^f\mbox{ verifies the Braid Equation}\}.$$
We then have
\begin{prop}
Let $(X,\trid)$ be a crossed set and $(C^n(X),\dt)$ be given as before. Then
\begin{enumerate}
\item $H^1(X)=\k^{\pi_0(X)}$, where $\pi_0(X)$ is the set of equivalence
classes of the relation generated by $j\sim i\trid j$.
\item $H^2(X)=\bdng(X)$ modulo change of basis of the form
$i\mapsto\lmb_i i$ for $\lmb_i\in\unid\k\ \forall i\in X$.
\end{enumerate}
\end{prop}
\begin{proof}
Straightforward (see \cite{gtesis}).
\end{proof}
Notice that we always have an inclusion $\unid\k\inclu H^2(X)$ by
$q\mapsto (f(i,j)=q)$.

\begin{lem}
Let $(X,\trid)$ be a crossed set and $f$ a $2$-cocycle in $C^2(X)$.
If $f$ takes values on the roots of unity then $(\k X,c)$ is of finite GT.
\end{lem}
\begin{proof}
We have to prove that the group-likes associated with the basis $X$ of
$\k X$ generate a finite group. These group-likes are nothing but
$g_i\in GL(\k X)$, $g_i(j)=f(i,j)i\trid j$. Let $N_{ij}=N(f(i,j))$, let
$N$ be the least common multiple of the $N_{ij}$'s and let $G\subset\unid\k$
be the group of $N$-roots of unity. Let $Y=X\times G$ and $\SS_Y$ the
symmetric group on the (finite) set $Y$. It is clear then that there
is an inclusion of the group generated by the $g_i$'s into $\SS_Y$.
\end{proof}

We can state \ref{rm:saq} in this setting: if $f$ is a $2$-cocycle and
$k\trid i=j$, then we have $f(i,i)=f(j,j)$. To see this, just compute
$$1=\dt(f)(k,i,i)=f(i,i)f(j,j)^{-1}f(k,i)^{-1}f(k,i).$$

\section{Rank $3$ Nichols algebras}\label{ss:r3}
In this section we show that if $V$ is a $3$-dimensional BP
of GT, then it is either of diagonal GT or
$V\simeq V_{3,q}$, $q$ a root of unity
(we give the definition of $V_{3,q}$ in \ref{df:v3q}).

Let $(V,c)$ be a $3$-dimensional BP of group type which is not of
abelian group type. Let $\{\bar x_0,\bar x_1,\bar x_2\}$ be a basis
as in \eqref{eq:bat} and $g_i\in\GL(V)$ the associated automorphisms
($i=0,1,2$). If two of the $g_i$'s were equal then $(V,c)$ would be
of abelian GT. For, suppose $g_1=g_2$. If $g_0=g_1$, the assertion
is obvious. If not, $g_0(\bar x_1)\in\spn\{\bar x_0\}$ or
$g_0(\bar x_1)\in\spn\{\bar x_1,\bar x_2\}$. In the first case,
we have $g_0g_1g_0^{-1}=g_0$ (which would imply that $g_0=g_1$),
while in the second case we have $g_0g_1g_0^{-1}=g_1$ (which would
imply that $g_0$ commutes with $g_1$ and hence the assertion).
Therefore, $\#\{g_0,g_1,g_2\}=3$. Consider now the subindices of the
$g_i$'s and the $x_i$'s to be in $\ZZ/3$. We have $g_ig_jg_i^{-1}=g_{-i-j}$,
for if not we would fall in one of the contradictions of before.
Furthermore, \eqref{eq:edc} implies that
$c(\bar x_i\otimes\bar x_j)=f(i,j)\bar x_{-i-j}\otimes\bar x_i$
for some $f:\{0,1,2\}\times\{0,1,2\}\to\unid\k$. That is, $(V,c)$ can
be constructed as in \ref{df:chcc}. Take $q=f(0,0)$ (i.e.,
$g_0(\bar x_0)=q\bar x_0$), and define $x_0=\bar x_0$,
$qx_1=g_2(x_0)$, $qx_2=g_0(x_1)$. Then it can be seen by hand
that $c(x_i\otimes x_j)=qx_{-i-j}\otimes x_i$. To see this, we must
compute $g_r(x_s)$ for $r,s\in\ZZ/3$, but this can be accomplished
just by using the definition of the $x_i$'s and the relations
$g_ig_jg_i^{-1}=g_{-i-j}$.
Alternatively, we can verify that $H^2(X)=\unid\k$ for $X=\ZZ/3$ and
$i\trid j=-i-j$.

\begin{defn}\label{df:v3q}
We denote this space by $V_{3,q}$. That is, $V_{3,q}$ has a basis
$\{x_0,x_1,x_2\}$ (the subindices in $\ZZ/3$) and the braiding is given by
$$c(x_i\otimes x_j)=qx_{-i-j}\otimes x_i.$$
Notice that $c^3=q^3$, from where the eigenvalues of $c$ belong to
$\{q,\xi q, \xi^2 q\}$ for $\xi$ a root of unity of order $3$.
\end{defn}

As a consequence, we have a classification of $3$-dimensional BP of
group type:
\begin{lem}
Let $(V,c)$ be a $3$-dimensional BP of GT. Then either
\begin{itemize}
\item $(V,c)$ is of abelian GT, or
\item $(V,c)=V_{3,q}$ for some $q\in\unid\k$.
\end{itemize}
Then, if $\Gm$ is a finite group and $V\in\ydskg$ is $3$-dimensional,
we have a basis $\{x_0,x_1,x_2\}$ such that either
\begin{itemize}
\item $c(x_i\otimes x_j)=q_{ij}x_j\otimes x_i$ for roots of unity $q_{ij}$, or
\item $c(x_i\otimes x_j)=qx_{-i-j}\otimes x_i$ (the subindices in $\ZZ/3$)
for $q$ a root of unity.
\end{itemize}

\end{lem}

The algebra $\toba{V_{3,-1}}$ has been studied in \cite{ms} and
is also included in \cite{ag}; it is $12$-dimensional and can
be presented with generators $\{x_0,x_1,x_2\}$ and relations
$$x_i^2=0\ (i=0,1,2),\quad
x_0x_1+x_1x_2+x_2x_0=0,\quad
x_0x_2+x_1x_0+x_2x_1=0.$$

Let $q$ be a root of unity of order $3$. By the remarks given above on the
eigenvalues of $c$, the morphism $1+c$ is injective, and hence
$\tobag 2{V_{3,q}}$ has dimension $9$. Since the image of $\ptl_0$ restricted
to $\tobag 2{V_{3,q}}$ has dimension $3$, its kernel has dimension $6$.
Furthermore, \ref{fc:ltdf} part \ref{fc:ltdf3} says that the set
$\{x_1^ix_2^j, 0\le i,j<3,\ i+j>2\}$ is linearly independent and it is
clearly contained in $\ker\ptl_0$. Taking into
account the vectors $1,x_1,x_2$ which lie also in $\ker\ptl_0$, we see
that the Hilbert polynomial of $\ker\ptl_0$ has coefficients greater or
equal to $P(t)=1+2t+6t^2+2t^3+t^4$, and then \ref{fc:ltdf} part \ref{fc:ltdf1}
says that $\dim\toba{V_{3,q}}\ge 3P(1)=36$.

\begin{rem}\label{rm:v3qsg}
By \ref{fc:ltdf} part \ref{fc:ltdf3}, $\dim(\toba{V_{3,q}})\ge N(q)^3$,
but the same result says that the equality does not hold. In fact,
$\dim(\toba{V_{3,-1}})>N(-1)^3$. Suppose then $q\neq -1$.
If we had $\dim(\toba{V_{3,q}})=N(q)^3$, then
\ref{fc:ltdf} part \ref{fc:ltdf3} would tell that
$\dim\tobag 2{V_{3,1}}=6$, or equivalently, that
$\dim\ker\SSS^2=3$. This is easily seen to be false, since for
$N(q)\neq 6$ we have $\ker\SSS^2=\ker(1+c)=\{0\}$, and for
$N(q)=6$ we have
$$\ker\SSS^2=\spn\{x_0x_1-qx_2x_0+q^2x_1x_2,\ 
x_1x_0-qx_2x_1+q^2x_0x_2\}.$$
\end{rem}

The unique Nichols algebra of rank $3$ which does not come
from the abelian case and whose dimension is $<32$ is hence
$\toba{V_{3,-1}}$. In particular, by \ref{fc:ltdf} part
\ref{fc:ltdf3} we see that it is the lowest dimensional
Nichols algebra which does not come from the abelian case.

\renewcommand{\theequation}{\arabic{section}.\arabic{subsection}.\arabic{equation}}
\section{Rank $4$ Nichols algebras}\label{ss:r4}
We classify now $4$-dimensional braided pairs of finite GT.
Let $V\in\ydskg$ be $4$-dimensional and suppose
$V=\oplus_iM_i=\oplus_iM(g_i,\rho_i)$ is its
decomposition as a sum of irreducible modules.

If $M(g,\rho)$ is $2$-dimensional, then it is of abelian GT.
If it is $1$-dimensional then $g$ is central in $\Gm$. Thus,
a sum of $1$-dimensional modules and at most one $2$-dimensional
module in $\ydskg$ is of abelian GT.
We have the following possibilities:
\begin{enumerate}
\item $4=1+1+1+1$, i.e. a sum of four $1$-dimensional
modules. It is of abelian GT.
\item $4=1+1+2$, it is also of abelian GT.
\item $4=1+3$. If the module of dimension $3$ is of
abelian GT, then $V$ is of abelian GT. The possibility
we must consider is that the $3$-dimensional module be
$V_{3,q}$ for some $q$. \label{r43m1}
\item $4=2+2$. Here $V$ is of abelian GT unless
$\rho_1$ and $\rho_2$ be characters and $g_1$
do not commute with $g_2$. \label{r42m2}
\item $4=4$. $V=M(g,\rho)$ is irreducible. There are
three sub-cases:
\begin{enumerate}
\item $g$ is central and $\deg\rho=4$. It is of abelian GT
(furthermore, by \ref{lm:csgr}, $\rho(g)=-1$ and $\dim\toba V=16$).
\item $[\Gm:\Gm_g]=2$ and $\deg\rho=2$. It is of abelian GT
(furthermore, by \ref{lm:csgr}, $N(\rho(g))=2$ or $3$.
If $N(\rho(g))=3$ we have by \ref{fc:ltdf} part \ref{fc:ltdf2}
that $\dim\toba V\ge 27^2$).
\item $[\Gm:\Gm_g]=4$ and $\rho$ is a character. \label{r41s}
\end{enumerate}
\end{enumerate}

We consider as before the situations which are not of abelian GT.

\subsection{Case \ref{r43m1}}\label{sss:r43m1}
Let $V=V_1\oplus V_{3,q}$ with $\dim V_1=1$, $V_1=M(g,\chi)$.
By \ref{fc:ltdf} part \ref{fc:ltdf2}, we have
$$\dim\toba V\ge\dim\toba{V_1}\dim\toba{V_{3,q}}\ge 36$$
unless $\chi(g)=q=-1$, case in which $\dim\toba V\ge 24$.
We denote by $\chi',g'$ the character and the element of
$\Gm$ giving $V_{3,-1}$, i.e. $V=M(g,\chi)\oplus M(g',\chi')$.
Since $g$ is central, it acts on $V_{3,-1}$ by a scalar,
say $q_{12}$. We denote $q_{21}=\chi(g')$. Let
$\sgm:=q_{12}q_{21}$. By \ref{fc:ltdf} part \ref{fc:ltdf2},
$\dim\toba V=24$ iff $\sgm=1$. If $\sgm\neq 1$, we consider
$x$ a generator of the space affording $\chi$ and
$x_0,x_1,x_2$ a basis of $M(g',\chi')$ as in section
\ref{ss:r3}, and define for $i=0,1,2$
$$z_i=\Ad_{x_i}(x)=x_ix-q_{21}xx_i.$$
Let $N=N(\sgm)$. It is straightforward (though tedious)
to see, using derivations, that the set
$$\{x^iz_0^{j_0}\ |\ 0\le i<2,\,0\le j_0<N\}$$
is linearly independent.
With this, we see also using derivations that the set
$$\{x^iz_0^{j_0}z_1^{j_1}\ |\ 0\le i<2,\,0\le j_0,j_1<N\}$$
is linearly independent, and finally we see that the set
$$\{x^iz_0^{j_0}z_1^{j_1}z_2^{j_2}\ |
	\ 0\le i<2,\,0\le j_0,j_1,j_2<N\}$$
is linearly independent. Furthermore, it is clear that
this set is contained in $K$, the intersection of the
kernels of $\ptl_{x_0},\ptl_{x_1},\ptl_{x_2}$. Thus, by
\ref{fc:ltdf} part \ref{fc:ltdf1} (taking $W=V_{3,q}$),
we see that $\dim\toba V\ge 2N^3\times 12\ge 192$.
Henceforth, in section \ref{ss:dm32} we shall not
consider these cases, except for $\sgm=1$.

\subsection{Case \ref{r42m2}}\label{sss:r42m2}
We have $V=M_1\oplus M_2=M(g_1,\rho_1)\oplus M(g_2,\rho_2)$,
$\rho_1,\rho_2$ are characters, $g_1,g_2$ do not commute, and
$[\Gm:\Gm_{g_1}]=[\Gm:\Gm_{g_2}]=2$. Let \linebreak
$g'_1=g_2g_1g_2^{-1}$,
$g'_2=g_1g_2g_1^{-1}$. The conjugacy class of $g_i$ is then
$\{g_i,g'_i\}$ and, by \cite[Lemma 3.1.9]{ag}, $[g_i,g'_i]=1$.
Let $H$ be the subgroup of $\Gm$ generated by $g_1,g'_1$,
which is commutative. Let $\sgm\in\Aut H$ be the restriction
of the adjoint $\sgm(h)=g_2hg_2^{-1}$.
We have $\sgm^2=\id$, and let $t=g'_1g_1^{-1}$. Notice that
\begin{align*}
g'_2&=g_1g_2g_1^{-1}=g_2g_2^{-1}g_1g_2g_1^{-1}=g_2g'_1g_1^{-1}=g_2t, \\
g_2&=g_1g'_2g_1^{-1}=g'_2g'_1g_1^{-1}=g'_2t=g_2t^2,
\end{align*}
whence $t^2=1$.
Let then $q_1=\rho_1(g_1)$, $q_2=\rho_2(g_2)$, $\eta_1=\rho_1(t)$,
$\eta_2=\rho_2(t)$ (thus $\eta_1=\pm 1$, $\eta_2=\pm 1$),
and $\al_1=\rho_1(g_2^2)$, $\al_2=\rho_2(g_1^2)$.
We have then the bases \linebreak
$\{x_1,x'_1=g_2\ganda x_1\}$ for $M_1$ and
$\{x_2,x'_2=g_1\ganda x_2\}$ for $M_2$, and the
braidings are given respectively by
$$\mbox{matrix of $c$ in }M_1
	=\mdpd {q_1}{\eta_1 q_1}{\eta_1 q_1}{q_1},\qquad
\mbox{matrix of $c$ in }M_2
	=\mdpd {q_2}{\eta_2 q_2}{\eta_2 q_2}{q_2},$$
and on the other monomials the braiding acts as
\begin{align*}
c(x_1\otimes x_2)=x'_2\otimes x_1, \qquad &
	c(x_2\otimes x_1)=x'_1\otimes x_2, \\
c(x_1\otimes x'_2)=\al_2 x_2\otimes x_1,\qquad  &
	c(x_2\otimes x'_1)=\al_1 x_1\otimes x_2, \\
c(x'_1\otimes x_2)=\eta_2x'_2\otimes x'_1,\qquad  &
	c(x'_2\otimes x_1)=\eta_1 x'_1\otimes x'_2, \\
c(x'_1\otimes x'_2)=\eta_2\al_2 x_2\otimes x'_1,\qquad  &
	c(x'_2\otimes x'_1)=\eta_1\al_1 x_1\otimes x'_2.
\end{align*}

We notice that, as in the proof of \ref{lm:csgr}, the
modules $M_i$ ($i=1,2$) are of Cartan type with associated
matrices $\mdpdc 2{2-N_i}{2-N_i}2$, where $N_i=N(q_i)$.
Thus $\toba V$ is infinite dimensional if $N_1>3$ or $N_2>3$.
Furthermore, if $N_1=3$ then $\toba{M_1}$ has dimension $27$
and analogously for $M_2$, whence, by \ref{fc:ltdf} part
\ref{fc:ltdf2}, if $N_1>2$ or $N_2>2$ then
$\dim\toba V\ge 4\times 27=108$.

We then consider $q_1=q_2=-1$. There are two cases:
$(\eta_1,\eta_2)=(1,1)$ and $(\eta_1,\eta_2)\neq (1,1)$.
In the case $(\eta_1,\eta_2)=(1,1)$, we define for
$\ep=\pm 1$ the vectors
$$z_{\ep}=\sqrt{\al_1}x_1+\ep x'_1,\quad
	z'_{\ep}=\sqrt{\al_2}x_2+\ep x'_2,$$
which give a basis of $V$. It is straightforward to compute the
braiding in this basis, which is given by
\begin{align*}
c(z_{\ep}\otimes z_{\ep'}) &= -z_{\ep'}\otimes z_{\ep}, \\
c(z_{\ep}\otimes z'_{\ep'}) &= \sqrt{\al_2}\ep' z'_{\ep'}\otimes z_{\ep}, \\
c(z'_{\ep'}\otimes z_{\ep}) &= \sqrt{\al_1}\ep z_{\ep}\otimes z'_{\ep'}, \\
c(z'_{\ep}\otimes z'_{\ep'}) &= -z'_{\ep'}\otimes z'_{\ep},
\end{align*}
from where we see that $\toba V$ comes from the abelian case. If we take
as a subspace $M_{\ep,\ep'}$ to be the linear span of $\{z_{\ep},z'_{\ep'}\}$
we have that the braiding has the matrix
$\mdpd{-1}{\sqrt{\al_2}\ep'}{\sqrt{\al_1}\ep}{-1}$.
By virtue of \ref{pr:mdr2}, $\toba{M_{\ep,\ep'}}$ has
dimension $4N(\ep\ep'\sqrt{\al_1}\sqrt{\al_2})$, from where
$$\dim(\toba V)\ge\dim(\toba{M_{\ep,\ep'}})\dim(\toba{M_{-\ep,-\ep'}})
	=16(N(\ep\ep'\sqrt{\al_1}\sqrt{\al_2}))^2\ge 64,$$
taking $\ep,\ep'$ in such a way that
$\ep\ep'\sqrt{\al_1}\sqrt{\al_2}\neq 1$.
Thus we will not be interested in these
algebras in section \ref{ss:dm32}.

\begin{rem}\label{rm:msvca}
A Nichols algebra as these can be found in \cite{ms},
where the authors consider a $4$-dimensional Yetter--Drinfeld
module over $\k\DD_4$ which can be decomposed as a sum of two
$2$-dimensional ones. That algebra is associated to the
constants $q_1=q_2=-1$, $\eta_1=\eta_2=1$, $\al_1=\al_2=1$.
It can be seen taking $M_{+,-}$ and $M_{-,+}$ that each
$\toba{M_{\pm,\mp}}$ is of type $A_2$ and hence
$8$-dimensional, and that for the elements of
$M_{\pm,\mp}\otimes M_{\mp,\pm}$ the braiding verifies
$c^2=\id$, from where, by \ref{fc:ltdf} part
\ref{fc:ltdf2}, $\dim\toba V=64$.
\end{rem}

\noindent If $(\eta_1,\eta_2)\neq(1,1)$, we suppose $\eta_1=-1$
(if $(\eta_1,\eta_2)=(1,-1)$ we interchange $M_1$ and $M_2$).
We take $z'_{\ep}$ as before, and define
\begin{align*}
t_+=\Ad_{x_1}z'_{+} &=x_1z'_{+}-\sqrt{\al_1}z'_{+}x_1, \\
t_-=\Ad_{x_1}z'_{-} &=x_1z'_{-}+\sqrt{\al_1}z'_{-}x_1.
\end{align*}
It is straightforward (but tedious) to see using derivations that
$t_{+},t_{-}$ are linearly independent, and that the set
$\{t_+^{i_+}t_-^{i_-}(z'_+)^{j_+}(z'_-)^{'j_-}\ |
	\ 0\le i_+,i_-,j_+,j_-<2\}$
is l.i. too. Since this set is contained in
$\ker\ptl_{x_1}\cap\ker\ptl_{x'_1}$, it is immediate by
a generalization of \ref{fc:ltdf} part \ref{fc:ltdf1} (see
\cite{gfr}) that the set
$$\{t_+^{i_+}t_-^{i_-}(z'_+)^{j_+}(z'_-)^{j_-}
	(x_1)^{t_1}(x'_1)^{t_2}
	\ |\ 0\le i_+,i_-,j_+,j_-,t_1,t_2<2\}$$
is linearly independent, whence $\dim\toba V\ge 64$. Hence,
in section \ref{ss:dm32} we shall not consider these algebras.

\subsection{Case \ref{r41s}}
We have $V=M(g,\rho)$, with $\orb_g=\{t_0=g,t_1,t_2,t_3\}$
and $\rho$ a character. We consider as in the
paragraph preceding \ref{rm:osf}, $f_g:\Gm\to\SS_4$.
By \ref{rm:osf} we have that the cardinality of
$f_g(\orb_g)$ may be $1,2,4$. If $\#f_g(\orb_g)=1$
then $\orb_g$ generates a commutative subgroup and
hence $V$ is of abelian GT.

If $\#f_g(\orb_g)=2$ we may index the $t_i$'s in such a way that
$$f_g(t_0)=f_g(t_1)=(2\ 3),\qquad f_g(t_2)=f_g(t_3)=(0\ 1)$$
(since $f_g(t_0)=f_g(t_1)$ then $t_0$ commutes with $t_1$, whence
either $f_g(t_0)=\id$ --which would imply $\#f_g(\orb_g)=1$--
or $f_g(t_0)=(2\ 3)$. Analogously for $f_g(t_2)$.).
Let $h\in\Gm$ be such that $ht_0h^{-1}=t_2$; then $f_g(h)(0)=2$.
It is easy to see (multiplying $h$ by the $t_i$'s on both sides
if necessary) that $\exists h'\in\Gm$ such that
$f_g(h')=(0\ 2)(1\ 3)$.
Let $x_0$ be a generator of the space affording $\rho$,
$x_1=t_2\ganda x_0$, $x_2=h'\ganda x_0$, $x_3=t_0h'\ganda x_0$.
It is straightforward to see that $\toba V$ is a particular
case of the algebras of \ref{sss:r42m2}. In fact, if
$q=\rho(t_0)$, $\eta=\rho(t_0^{-1}t_1)$ and $\al=\rho(t_2^2)$
then $\toba V\simeq\toba{M_{q,q,\eta,\eta,\al,\al}}$,
and then $\dim\toba V\ge 64$. We shall not consider them
in section \ref{ss:dm32}.

\bigskip

If $\#f_g(\orb_g)=4$, then $f_g(t_i)$ fixes $i$ and does not fix $j$ for
$j\neq i$, and then it is a tri-cycle. We may index the $t_i$'s in such a
way that $f_g(t_0)=(1\ 2\ 3)$. It can be seen that the hypothesis
$f_g(t_1)=(0\ 2\ 3)$ leads to a contradiction,
and thence we have $f_g(t_1)=(0\ 3\ 2)$.
This implies $f_g(t_2)=(0\ 1\ 3)$ and $f_g(t_3)=(0\ 2\ 1)$.
Take $h=t_2t_1t_0^{-2}$. Then $f_g(h)=\id$, and thus $h$ commutes with
$\orb_g$ (notice that for the same reason $t_i^3$ commutes with $\orb_g$
for $i=0,1,2,3$). We note that $t_1^3=t_2t_0^3t_2^{-1}=t_0^3$, and
in the same way we see that $t_i^3=t_j^3$. Notice also that
\begin{align*}
h&=t_2t_1t_0^{-2}=t_0^{-2}t_0^2(t_2t_1)t_0^{-2}
	=t_0^{-2}t_1t_3,\mbox{ and} \\
h^2&=t_0^{-2}t_1t_3t_2t_1t_0^{-2}=t_0^{-2}t_1t_3t_2t_3^{-2}t_1
	=t_0^{-2}t_1t_1t_3t_3^{-2}t_1 \\
&=t_0^{-2}t_1^2t_3^{-1}t_1=t_0^{-2}t_0^{-1}t_1^2t_1=t_0^{-3}t_1^3=1.
\end{align*}
Take $x_0$ a generator of the space affording $\rho$,
let $q=\rho(g)$ and take $qx_1=t_2\ganda x_0$, $qx_2=t_0\ganda x_1$,
$qx_3=t_0\ganda x_2$. Let $\al=\rho(h)=\pm 1$; we compute
\begin{align*}
t_0\ganda x_0 &= qx_0,\quad t_0\ganda x_1=qx_2,\quad
t_0\ganda x_2=qx_3, \\
t_0\ganda x_3 &= q^{-2}t_0^3\ganda x_1=q^{-3}t_0^3t_2\ganda x_0
	=q^{-3}t_2t_0^3\ganda x_0=t_2\ganda x_0=qx_1,  \\
t_1\ganda x_0 &= t_0^2t_2t_0^{-2}\ganda x_0=q^{-2}t_0^2t_2\ganda x_0
	=q^{-1}t_0^2\ganda x_1=qx_3, \\
t_1\ganda x_1 &= qx_1 \mbox{ by \ref{rm:saq}}, \\
t_1\ganda x_2 &= t_2^{-1}t_2t_1t_0^{-2}t_0^2\ganda x_2=q^2t_2^{-1}h\ganda x_1
	=q^2\al t_2^{-1}\ganda x_1=q\al x_0, \\
t_1\ganda x_3 &= t_0^2t_0^{-2}t_1t_3t_3^{-1}\ganda x_3=q^{-1}\al t_0^2\ganda x_3
	=q\al x_2, \\
t_2\ganda x_0 &= qx_1, \\
t_2\ganda x_1 &= t_2t_1t_0^{-2}t_0^2t_1^{-1}\ganda x_1
	=q^{-1}t_2t_1t_0^{-2}t_0^2\ganda x_1=q\al x_3, \\
t_2\ganda x_2 &= qx_2 \mbox{ by \ref{rm:saq}}, \\
t_2\ganda x_3 &= t_2^3t_2^{-2}\ganda x_3=q^{-1}\al^{-1} t_2^3t_2^{-1}\ganda x_1
	=q^{-2}\al^{-1}t_2^3\ganda x_0 \\
	&=q^{-2}\al t_0^3\ganda x_0=q\al x_0,
\end{align*}
\begin{align*}
t_3\ganda x_0 &= t_0t_2t_0^{-1}\ganda x_0=t_0\ganda x_1=qx_2, \\
t_3\ganda x_1 &= t_0t_2t_0^{-1}\ganda x_1=q^{-1}t_0t_2\ganda x_3
	=\al t_0\ganda x_0=q\al x_0, \\
t_3\ganda x_2 &= t_0t_2t_0^{-1}\ganda x_2=q^{-1}t_0t_2\ganda x_1
	=\al t_0\ganda x_3=q\al x_1, \\
t_3\ganda x_3 &= qx_3 \mbox{ by \ref{rm:saq}}.
\end{align*}

We shall denote by $V_{4,q,\al}$ the $4$-dimensional module generated by
$\{x_0,x_1,x_2,x_3\}$ with the braiding of above. By \ref{fc:ltdf}
part \ref{fc:ltdf3}, we have $\dim(\toba{V_{4,q,\al}})\ge 81$ if
$N(q)\ge 3$, whence we will be interested in the case $q=-1$.

\begin{algo}\label{al:cdt72}
If $\al=-1$, it is easy to see that for $V_{4,-1,-1}$ the restriction
of $c$ to the linear span of $\{x_i\otimes x_j\ |\ i\neq j\}$ verifies
$c^3=1$, from where
$$\ker\SSS^2=\ker(1+c)=\spn\{x_i^2,\ i=0,1,2,3\}.$$
This implies that $\dim\tobag 2{V_{4,-1,-1}}=12$. Now, since the image
of the map \linebreak
$\ptl_0:\tobag 2{V_{4,-1,-1}}\to\tobag 1{V_{4,-1,-1}}$ is
$3$-dimensional (we have $x_0\notin\im(\ptl_0)$ because $x_0^2=0$),
we have $\dim(\ker(\ptl_0)\cap\tobag 2{V_{4,-1,-1}})=9$. Thus, the
Hilbert polynomial of $\toba {V_{4,-1,-1}}$ verifies
$$P_{\toba {V_{4,-1,-1}}}(t)\ge (1+t)(1+3t+9t^2+3t^3+t^4)$$
by virtue of \ref{fc:ltdf2}, whence $\dim\toba {V_{4,-1,-1}}\ge 34$.
We shall not consider this algebra in section \ref{ss:dm32}.

For the module $V_{4,-1,1}$, it is easy to find the relations in
degree $2$:
\begin{align*}
&x_i^2=0,\ i=0,1,2,3 \\
&x_0x_1+x_2x_0+x_1x_2=0, \\
&x_0x_2+x_3x_0+x_2x_3=0, \\
&x_0x_3+x_1x_0+x_3x_1=0, \\
&x_1x_3+x_2x_1+x_3x_2=0.
\end{align*}
With the order $x_0>x_1>x_2>x_3$, to give a Gr\"obner basis of the ideal
generated we have to add the relations
\begin{align*}
x_1x_2x_1-x_2x_1x_2 &=0, \\
x_2x_3x_2-x_3x_2x_3 &=0.
\end{align*}
It is not so easy\footnote{I used a computer program by
E. M\"uller \cite{mul} to see that there should be one more relation
in degree $6$.} to see that there is one more relation for
$\toba{V_{4,-1,1}}$:
$$x_2x_3x_4x_2x_3x_4+x_3x_4x_2x_3x_4x_2+x_4x_2x_3x_4x_2x_3=0.$$
With this, it is straightforward (though a little bit tedious)
to see that the dimension of $\tobag i{V_{4,-1,1}}$ is
$1,4,8,11,12,12,11,8,4,1$ for $i=0,1,2,3,4,5,6,7,8,9$ respectively,
and hence $\dim(\toba{V_{4,-1,1}})=72$.
We thus shall not consider this algebra.
\end{algo}

\renewcommand{\theequation}{\arabic{section}.\arabic{equation}}
\section{Classification of Nichols algebras of dimension $<32$}\label{ss:dm32}

The Nichols algebras (and also their liftings) of dimension $p$ and
$p^2$ are classified in \cite{gfr}.
We shall consider now the other possibilities.

\medskip
Let $(V,c)$ be a BP of finite GT with Nichols rank $p^3$.
We have $r(n)=3$ (see \ref{df:rd} for the definition of $r(n)$),
from where $\dim V\le 3$. By the results in section \ref{ss:r3},
$V$ is of abelian GT, for, if not, we would have $V=V_{3,q}$,
but then by \ref{fc:ltdf} part \ref{fc:ltdf2}, $N(q)=p$ and by
\ref{rm:v3qsg}, the Nichols rank of $V$ would be $>p^3$,
a contradiction. Hence, if $\dim V=3$, $\toba V$ is
a QLS because of \ref{fc:ltdf} part \ref{fc:ltdf2}, i.e., $V$ has
a basis $\{x_1,x_2,x_3\}$ such that the set
$\{x_1^{n_1}x_2^{n_2}x_3^{n_3}\ |\ 0\le n_i<p\}$ is a basis of
$\toba V$. If $\dim V=2$ then $c$ has a matrix
$\mdpd{q_{11}}{q_{12}}{q_{21}}{q_{22}}$ where
the $q_{ii}$'s are either $p$-roots or $p^2$-roots of unity.
If one of them has order $p^2$, say $q_{22}$, then $N(q_{11})=p$
and $\toba V$ is a QLS. If $N(q_{11})=N(q_{22})=p$ then by
\ref{pr:mdr2} we have two possibilities:
\begin{enumerate}
\item $p=2$ and $q_{12}q_{21}\neq 1$. Furthermore,
for $\toba V$ to be $8$-dimensional we must have
$q_{12}q_{21}=-1$. The algebra $\toba V$ is of type $A_2$.
\item $p>2$ and $q_{12}q_{21}q_{22}=1$. Interchanging
$q_{11}$ with $q_{22}$ we see that also \linebreak
$q_{12}q_{21}q_{11}=1$, whence $q_{11}=q_{22}$ and
$\toba V$ is also of type $A_2$.
\end{enumerate}
Thus we reach in both cases to the same situation.
The algebra $\toba V$ has a PBW basis given by
$\{x_1^{m_1}z^nx_2^{m_2}\ |\ 0\le m_1,m_2,n<p\}$,
$z=x_2x_1-q_{21}x_1x_2$. The last case to be considered
is $\dim V=1$. We have simply $\{x\}$ a basis of $V$
such that $c(x\otimes x)=qx\otimes x$, $N(q)=p^3$.
Thus $\toba V=\k[x]/(x^{p^3})$.

We then have proved:
\begin{thm}
Let $(V,c)$ be a BP of finite group type with Nichols rank $p^3$
($p$ a prime number). Then $V$ is of diagonal GT and it fits
in one of the following possibilities:
$$\begin{array}{|r|l|} \hline
\dim V & \mbox{Type}  \\ \hline
3 & \mbox{QLS, }N(q_{11})=N(q_{22})=N(q_{33})=p \\ \hline
2 & \mbox{QLS, }N(q_{11})=p,\ N(q_{22})=p^2 \\ \hline
2 & A_2,\ q_{11}=q_{22}=(q_{12}q_{21})^{-1},\ N(q_{11})=p \\ \hline
1 & \mbox{QLS, }N(q_{11})=p^3 \\ \hline
\end{array}$$
\qed
\end{thm}

\medskip
Let $n=p_1p_2$ a product of two distinct primes such that
$p_1<p_2<p_1^2$. Then $r(n)<3$, whence $\dim V\le 2$ and
it is of abelian GT. If $\dim V=2$ and $c$ has matrix
$\mdpd{q_{11}}{q_{12}}{q_{21}}{q_{22}}$, then by \ref{fc:ltdf}
part \ref{fc:ltdf2} we must have (reversing the basis if necessary)
$N(q_{11})=p_1$, $N(q_{22})=p_2$, and $q_{12}q_{21}=1$,
i.e., $\toba V$ is a QLS with PBW
basis $\{x_1^{m_1}x_2^{m_2}\ |\ 0\le m_i<p_i,\ i=1,2\}$. If $\dim V=1$
then as before $x$ generates $V$, $c(x\otimes x)=qx\otimes x$ with
$N(q)=p_1p_2$ and $\toba V=\k[x]/(x^{p_1p_2})$.

Since we have computed the dimension of $\toba {V_{3,-1}}$,
together with the fact that we have shown that
$\dim\toba{V_{3,q}}\ge 32$ if $q\neq -1$ and $\dim\toba V\ge 32$
for $V$ a $4$-dimensional BP of finite GT which is not of abelian
GT, this last argument applies also for $p_1=2$ and $p_2<16$
(i.e., $n=10$, $n=14$, $n=22$ and $n=26$). In these cases,
if $\dim V$ were $3$, then $V$ would be of abelian GT because
if not $V=V_{3,-1}$, which would give a Nichols algebra of
dimension $12$, a contradiction. Then remark \ref{rm:csm}
applies and $\dim V\le 2$.

We have hence the classification of the Nichols algebras of dimensions
$1,2,3$, $4,5,6,7,8,9,10,11,13,14,15,17,19,21,22,23,25,26,27,29,31$.
For the remaining  \linebreak
dimensions it is routine to make the computations, and we refer
to \cite{gtesis}. The result can be stated as follows:

\medskip
$$\begin{array}{|r|c|} \hline
\mbox{Type} & \mbox{Dimensions}  \\ \hline
\mbox{QLS, rank $1$} &  \mbox{All dimensions} \\ \hline
\mbox{QLS, rank $2$} &  \neq p \\ \hline
\mbox{QLS,\ rank $3$} &  8,\ 12,\ 16,\ 18,\ 20,\ 24,\ 27,\ 28,\ 30 \\ \hline
\mbox{QLS,\ rank $4$} &  16,\ 24 \\ \hline
\mbox{$A_2$} &  8,\ 12,\ 16,\ 20,\ 27,\ 28 \\ \hline
\mbox{$A_2\times A_1$} & 16,\ 24 \\ \hline
\mbox{$V_{3,-1}$} & 12 \\ \hline
\mbox{$V_{3,-1}\times A_1$} & 24 \\ \hline
\end{array}$$

\medskip
When applying the lifting procedure, the computations in the previous
section are the first step in order to classify pointed Hopf algebras
of index $<32$. The following lemmas, which can be proved as in
\cite[Lemma 7.2]{as4} or \cite[\S 7]{g32} and are referred to
\cite{gtesis}, show that any such Hopf algebra is generated by
group-likes and skew-primitive elements, completing the second step
of the procedure.

\begin{lem}\label{lm:en1}
Let $(V=(x_1,x_2),c)$ be an BP of abelian GT
with matrix $(q_{ij})$, and suppose $\dim\toba V<32$.
If $S$ is a graded braided Hopf algebra between $TV$ and $\toba V$, i.e.
$$TV\sobre S\sobre\toba V,$$
with graded maps and categorical maps (that is, there exists a braided
category such that $TV,\ S$ and $\toba V$ belong to it and the projections
are maps in the category), then either $S=\toba V$ or $\dim S\ge 32$.
\qed\end{lem}

\begin{lem}
Let $S$ be a finite dimensional graded braided Hopf algebra with graded maps
$T(V_{3,-1})\sobre S\displaystyle\mathop{\sobre}^p\toba{V_{3,-1}}$.
Then $p$ is an isomorphism.\qed
\end{lem}

\end{document}